\documentclass[twocolumn]{svjour3} 

\smartqed  
\usepackage{graphicx}
\usepackage{amssymb,amsmath,amsfonts,latexsym,psfrag,mathptmx,url}
\usepackage{amsmath}

\newtheorem{algorithm}{Algorithm}

\newcommand{\scal}[2]{{#1}^\top{#2}}

\begin{document}

\title{Smart depth of field optimization applied to a robotised view camera}

\author{St\'ephane Mottelet \and Luc de Saint Germain \and Olivier Mondin}


\institute{S. Mottelet \at
               Laboratoire de Math\'ematiques Appliqu\'ees\\
Universit\' e de Technologie de Compi\`egne\\
60205 Compi\`egne France\\
              \email{stephane.mottelet@utc.fr}           
           \and
           L. de Saint Germain, O. Mondin\at
           Luxilon\\
21, rue du Calvaire\\
92210 Saint-Cloud France\\
\email{lsg@lsg-studio.com}
}
\maketitle

\begin{abstract}
The great flexibility of a view camera allows the acquisition of high quality images that would not be possible
any other way. Bringing a given object into focus is however a long and tedious task, although the underlying optical laws are known. 
A fundamental parameter is the aperture of the lens entrance pupil because it directly affects the depth of field. The smaller the aperture, the larger the depth of field. However a too small aperture destroys the sharpness of the image because of diffraction on the pupil edges. Hence, the desired optimal configuration of the camera is such that the object has the sharpest image with the greatest possible lens aperture. In this paper, we show that when the object is a convex polyhedron, an elegant solution to this problem can be found. The problem takes the form of a constrained optimization problem, for which theoretical and numerical results are given. 

\keywords{Large format photography \and Computational photography \and Scheimpflug principle}
\end{abstract}

\section{Introduction}
\begin{figure}
\begin{center}
\begin{tabular}{c}
\includegraphics[width=0.7\linewidth]{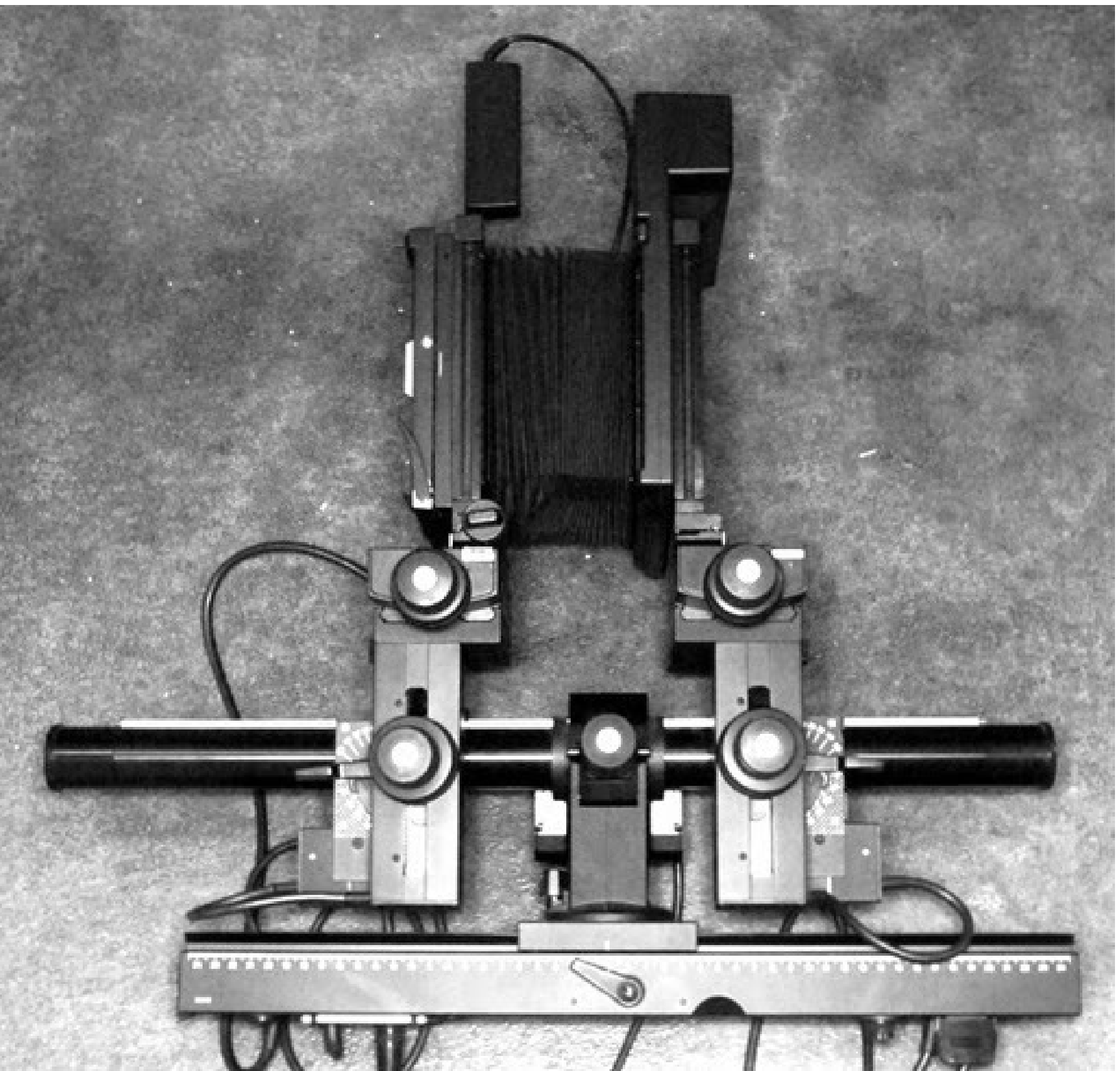}\\(a)\\\includegraphics[width=0.7\linewidth]{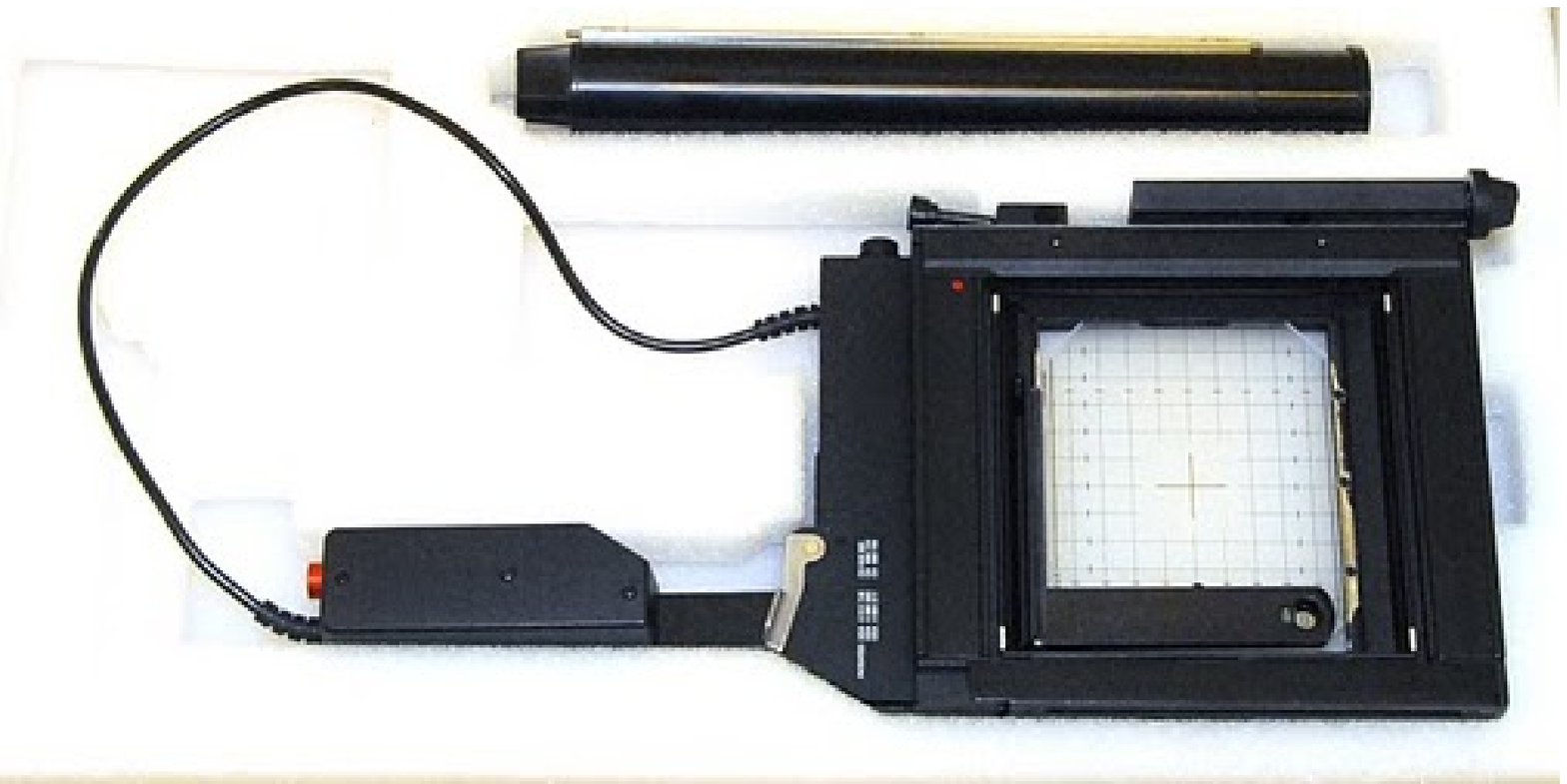}\\(b)
\end{tabular}
\end{center}
   \caption{The Sinar e (a) and its metering back (b).}
\label{fig:sinar_e}
\end{figure}

Since the registration of Theodor Scheimpflug's patent in 1904 (see \cite{Scheimpflug}), and the book of Larmore in 1965 where a proof of the so-called {\em Scheimpflug principle} can be found (see \cite[p. 171-173]{larmore}), very little has been written about the mathematical concepts used in modern view cameras, until the development of the Sinar e in 1988 (see Figure \ref{fig:sinar_e}). A short description of this camera is given in \cite[p. 23]{Tillmans}:

\medskip{\em The {\rm Sinar e} features an integrated electronic computer, and in the studio offers a maximum of convenience and optimum computerized image setting. The user-friendly software guides the photographer through the shot without technical confusion. The photographer selects the perspective (camera viewpoint) and the lens, and chooses the areas in the subject that are to be shown sharp with a probe. From these scattered points the Sinar e calculates the optimum position of the plane of focus, the working aperture needed, and informs the photographer of the settings needed}

\medskip Sinar sold a few models of this camera and discontinued its development in the early nineties. Surprisingly, there has been very little published user feedback about the camera itself. However many authors started to study (in fact, re-discover) the underlying mathematics (see e.g. \cite{Merklinger} and the references therein). The most crucial aspect is the consideration of depth of field and the mathematical aspects of this precise point are now well understood. When the geometrical configuration of the view camera is precisely known, then the depth of field region (the region of space where objects have a sharp image) can be determined by using the laws of geometric optics. Unfortunately, these laws can only be used as a rule of thumb
when operating by hand on a classical view camera.  Moreover, the photographer is rather interested in the inverse problem: given an object which has to be rendered sharply, what is the optimal configuration of the view camera? A fundamental parameter of this configuration is the aperture of the camera lens. Decreasing the lens aperture diameter increases the depth of field but also increases the diffraction of light by the lens entrance pupil. Since diffraction decreases the sharpness of the image, the optimal configuration should be such that the object fits the depth of field region with the greatest aperture. 

This paper presents the mathematical tools used in the software of a computer controlled view camera solving this problem. Thanks to the high precision machining of its components, and to the known optical parameters of the lens and digital sensor, a reliable mathematical model of the view camera has been developed. This model allows the acquisition of 3D coordinates of the object to be photographed, as explained in Section \ref{basics}. 
In Section \ref{focus} we study the depth of field optimization problem from a theoretical and numerical point of view. We conclude and briefly describe the architecture of the software in Section \ref{conclusion}.

\section{Basic mathematical modeling}
\label{basics}
\subsection{Geometrical model of the View Camera}
\begin{figure}[tp]
\begin{center}
\begin{tabular}{c}
\includegraphics[width=0.8\linewidth]{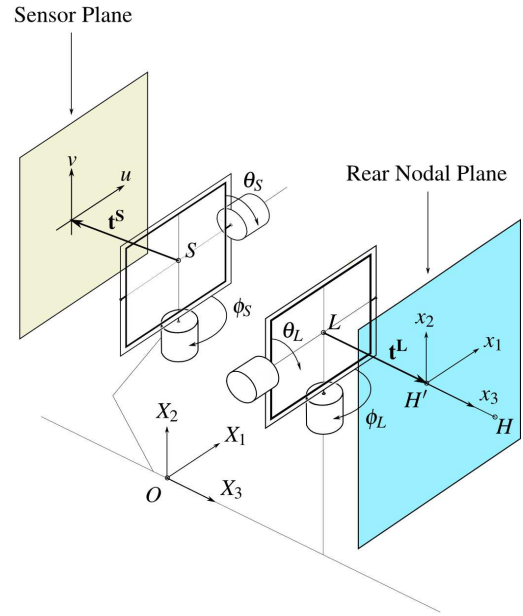}\\
(a)\\\\
\includegraphics[width=0.8\linewidth]{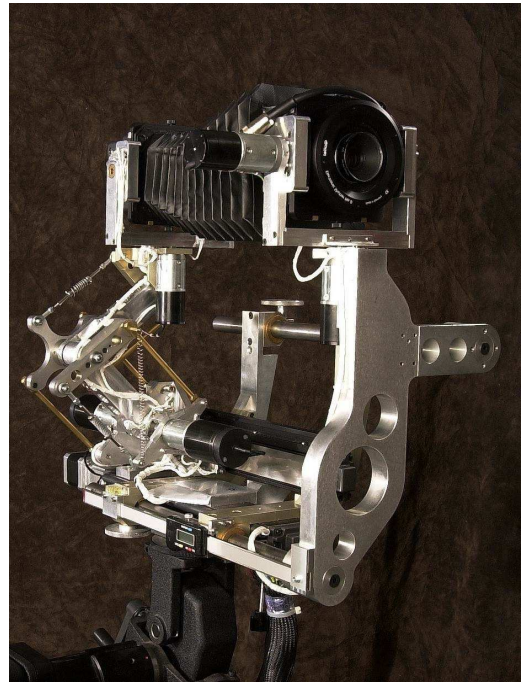}\\
(b)
\end{tabular}
\end{center}
   \caption{Geometrical model (a) and robotised view camera (b).}
\label{fig:viewcamera}
\end{figure}
We consider the robotised view camera depicted in Figure \ref{fig:viewcamera}(a) and its geometrical model in Figure \ref{fig:viewcamera}(b). We use a global Euclidean coordinate system $(\mathbf{O},X_1,X_2,X_3)$ attached to the camera's tripod. The front standard, symbolized by its frame with center $L$ of global coordinates $\mathbf{L}=(L_1,L_2,L_3)^\top$, can rotate along its tilt and swing axes with angles $\theta_L$ and $\phi_L$. Most camera lenses are in fact thick lenses and nodal points $H'$ and $H$ have to be considered (see \cite{Ray} p. 43-46).  The rear nodal plane, which is parallel and rigidly fixed to the front standard, passes through the rear nodal point $H'$. Since $L$ and $H'$ do not necessarily coincide, the translation between these two points is denoted $\mathbf{t^L}$. The vector $\overrightarrow{HH'}$ is supposed to be orthogonal to the rear nodal plane.

The rear standard is symbolized by its frame with center $S$, whose global coordinates are given by $\mathbf{S}=(S_1,S_2,S_3)^\top$. It can rotate along its tilt and swing axes with angles $\theta_S$ and $\phi_S$. The sensor plane is parallel and rigidly fixed to the rear standard. The eventual translation between $S$ and the center of the sensor is denoted by $\mathbf{t^S}$.

The rear standard center $S$ can move in the three $X_1,X_2$ and $X_3$ directions but the front standard center $L$ is fixed. The rotation matrices associated with the front and rear standard alt-azimuth mounts are respectively given by $\mathbf{R^L}=\mathbf{R}(\theta_L,\phi_L)$ and $\mathbf{R^S}=\mathbf{R}(\theta_S,\phi_S)$ where 
 $$
 \mathbf{R}(\theta,\phi)=\left(\begin{array}{ccc}\cos\phi & -\sin\phi\sin\theta & -\sin\phi\cos\theta\\
 0 &\cos\theta & -\sin\theta\\
 -\sin\phi&\cos\phi\sin\theta&\cos\phi\cos\theta\end{array}\right).
 $$
The intrinsic parameters of the camera (focal length $f$, positions of the nodal points $H$, $H'$,  translations $\mathbf{t^S}$, $\mathbf{t^L}$, image sensor characteristics) are given by their respective manufacturers data-sheets. The extrinsic parameters of the camera are $\mathbf{S}$, $\mathbf{L}$, the global coordinate vectors of $S$ and $L$, and the four rotation angles $\theta_S$, $\phi_S$, $\theta_L$, $\phi_L$. The precise knowledge of the extrinsic parameters is possible thanks to the computer-aided design model used for manufacturing the camera components. In addition, translations and rotations of the rear and front standards are controlled by stepper motors whose positions can be precisely known. In the following, we will see that this precise geometrical model of the view camera allows one to solve various photographic problems. The first problem is the determination of coordinates of selected points of the object to be photographed.

In the sequel, for sake of simplicity, we will give all algebraic details of the computations for a thin lens, i.e. when the two nodal points $H$, $H'$ coincide. In this case, the nodal planes are coincident in a so-called \textit{lens plane}. We will also consider that $\mathbf{t^L}=(0,0,0)^\top$ so that $L$ is the optical center of the lens. Finally, we also consider that $\mathbf{t}^S=(0,0,0)^\top$ so that $S$ coincides with the center of the sensor.

\subsection{Acquisition of object points coordinates}
\begin{figure}
\begin{center}
\begin{tabular}{c}
\includegraphics[width=\linewidth]{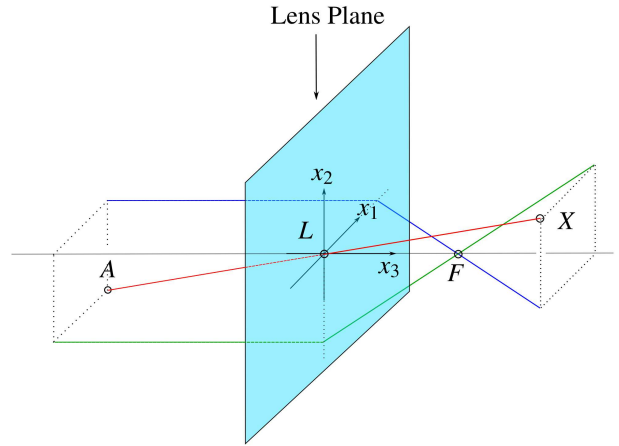}
\end{tabular}
\end{center}
\caption{Graphical construction of the image $A$ of an object point $X$ in the optical coordinate system. Green rays and blue rays respectively lie in the $(L,x_3,x_2)$ and $(L,x_3,x_1)$ planes and $F$ is the focal point.}
\label{lensmodel4}
\end{figure}
\begin{figure}
\begin{center}
\begin{tabular}{c}
\includegraphics[width=\linewidth]{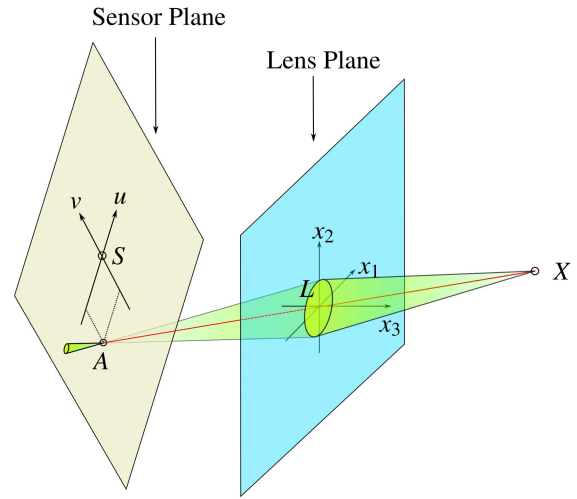}
\end{tabular}
\end{center}
\caption{Image formation when considering a lens with a pupil.}
\label{lensmodel5}
\end{figure}
Let us consider a point $X$ with global coordinates given by $\mathbf{X}=(X_1,X_2,X_3)^\top$. The geometrical construction of the image $A$ of $X$ through the lens is depicted in Figure \ref{lensmodel4}. We have considered a local optical coordinate system attached to the lens plane with origin $L$. The local coordinates of $X$ are given by $\mathbf{x}=(\mathbf{R^L})^{-1}(\mathbf{X}-\mathbf{L})$ and the focal point $F$ has local coordinates $(0,0,f)^\top$. Elementary geometrical optics (see \cite{Ray} p. 35-42) allows one to conclude that
if the local coordinates of $A$ are given by $\mathbf{a}=(a_1,a_2,a_3)^\top$, then $a_3$, $x_3$ and $f$ are linked by the thin lens equation given in its Gaussian form by 
$$
-\frac{1}{a_3}+\frac{1}{x_3}=\frac{1}{f}.
$$
Since $A$ lies on the $(XL)$ line, the other coordinates are obtained by straightforward computations and we have the conjugate formulas
\begin{align}
\mathbf{a}&=\frac{f}{f-x_3}\mathbf{x},\\
\mathbf{x}&=\frac{f}{f+a_3}\mathbf{a}.
\label{eq:conj}
\end{align}
Bringing an object into focus is one of the main tasks of a photographer but it can also be used to calculate the coordinates of an object point. It is important to remember that all light rays emanating from $X$ converge to $A$ but pass through a pupil (or diaphragm) assumed to be circular, as depicted in Figure \ref{lensmodel5}.
Since all rays lie within the oblique circular cone of vertex $A$ and whose base is the pupil, the image of $X$ on the sensor will be in focus only if the sensor plane passes through $A$, otherwise its extended image will be a blur spot. By using the full aperture of the lens, the image will rapidly go out of focus if the sensor plane is not correctly placed, e.g. by translating $S$ into the $x_3$ direction. This is why auto-focus systems on classical cameras only work at near full aperture: the distance to an object is better determined when the depth of field is minimal.  

The uncertainty on the position of $S$ giving the best focus is related to the diameter of the so-called ``circle of confusion", i.e. the maximum diameter of a blur spot that is indistinguishable from a point. Hence, everything depends on the size of photosites on the sensor and on the precision of the focusing system (either manual or automatic). This uncertainty is acceptable and should be negligible compared to the uncertainty of intrinsic and extrinsic camera parameters.

The previous analysis shows that the global coordinates of $X$ can be computed, given the position $(u,v)^\top$ of its image $A$ on the sensor plane. This idea has been already used on the Sinar e, where the acquisition of $(u,v)^\top$ was done by using a mechanical metering unit (see Figure \ref{fig:sinar_e} (b)). In the system we have developed, a mouse click in the live video window of the sensor is enough to indicate these coordinates. Once $(u,v)^\top$ is known, the coordinates of $A$ in the global coordinate system are given by
$$
\mathbf{A}=\mathbf{S}+\mathbf{R^S}\left(\begin{array}{c}u\\v\\0\end{array}\right),
$$
and its coordinates in the  optical system by
$$
\mathbf{a}=(\mathbf{R^L})^{-1}(\mathbf{A}-\mathbf{L}).
$$
Then the local coordinate vector $\mathbf{x}$ of the reciprocal image is computed with (\ref{eq:conj}) and the global coordinate vector $\mathbf{X}$ is obtained by
$$
\mathbf{X}=\mathbf{L}+\mathbf{R^L}\mathbf{x}.
$$
By iteratively focusing on different parts of the object, the photographer can obtain a set of points $\mathcal{X}=\{{X}^1,\dots,{X}^n\}$, with $n\ge 3$,  which can be used to determine the best configuration of the view camera, i.e. the positions of front and rear standards and their two rotations, in order to satisfy focus requirements.

\section{Focus and depth of field optimization}
\label{focus}
In classical digital single-lens reflex (DLSR) cameras, the sensor plane is always parallel to the lens plane and to the plane of focus. For example, bringing into focus a long and flat object which is not parallel to the sensor needs to decrease the aperture of the lens in order to extend the depth of field. On the contrary, view cameras with tilts and swings (or DLSR with a tilt/shift lens) allow to skew away the plane of focus from the parallel in any direction. Hence, bringing into focus the same long and flat object with a view camera can be done at full aperture. This focusing process is unfortunately very tedious. However, if a geometric model of the camera and the object are available, the adequate rotations can be estimated precisely. In the next sections, we will explain how to compute the rear standard position and the tilt and swing angles of both standards to solve two different problems: 
\begin{enumerate}
\item when the focus zone is roughly flat, and depth of field is not a critical issue, then the object plane is computed from the set of object points $\mathcal{X}$. If $n=3$ and the points are not aligned then this plane is uniquely defined. If $n>3$ and at least $3$ points are not aligned, we compute the best fitting plane minimizing the sum of squared orthogonal distances to points of $\mathcal{X}$. Then, we are able to bring this plane into sharp focus by acting on:
\begin{enumerate}
\item the angles $\theta_L$ and $\phi_L$ of the front standard and the position of the rear standard, for arbitrary rotation angles $\theta_S$, $\phi_S$.\label{cas1}
\item the angles $\theta_S$, $\phi_S$ and position of the rear standard, for arbitrary rotation angles $\theta_L$, $\phi_L$ (in this case there is a perspective distortion).
\end{enumerate}
\item when the focus zone is not flat, then the tridimensional shape of the object has to be taken into account. 
\end{enumerate}
The computations in case \ref{cas1} are detailed in Section \ref{section4}. In Section \ref{section5} a general algorithm is described that allows the computation of angles $\theta_L$ and $\phi_L$ of the front standard and the position of the rear standard such that all the object points are in the depth of field region with a maximum aperture. We give a theoretical result showing that the determination of the solution amounts to enumerate a finite number of configurations.

\subsection{Placement of the plane of sharp focus by using tilt and swing angles}
\label{section4}
\begin{figure}
\begin{center}
\includegraphics[width=\linewidth]{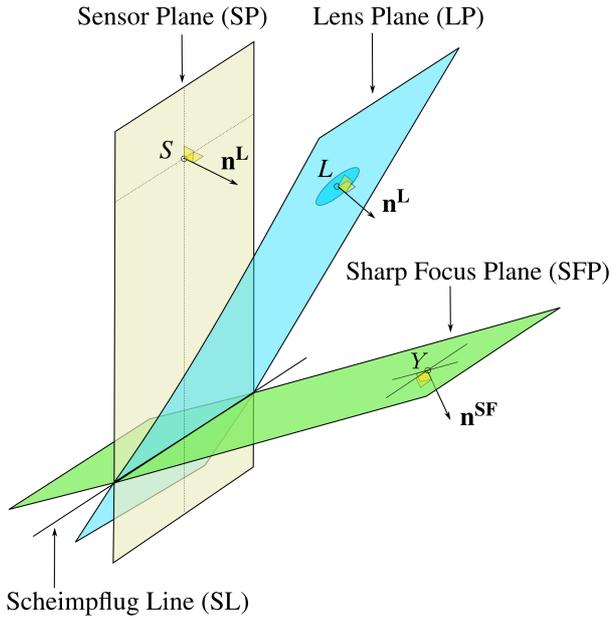}\\
\end{center}
\caption{Illustration of the Scheimpflug rule.}
\label{scheimpflug}
\end{figure}

\begin{figure}
\begin{center}
\includegraphics[width=\linewidth]{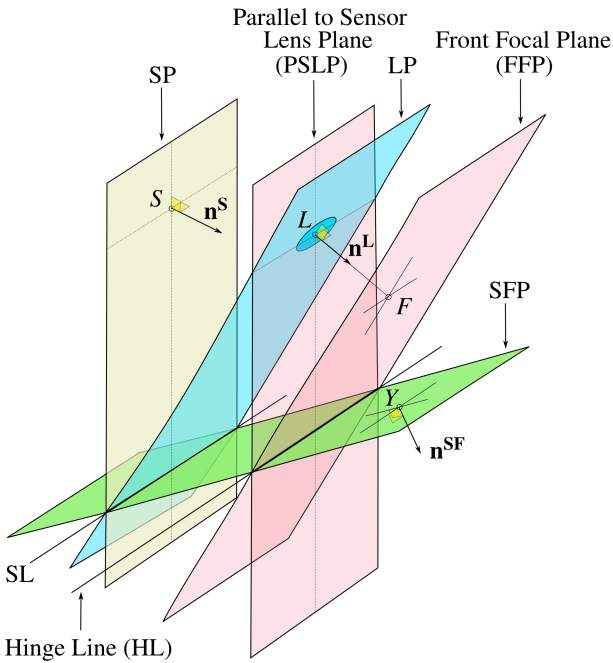}\\
\end{center}
\caption{Illustration of the Hinge rule.}
\label{hinge}
\end{figure}

In this section we study the problem of computing the tilt and swing angles of front standard and the position of the rear standard for a given sharp focus plane. Although the underlying laws are well-known and are widely described (see \cite{Merklinger,Wheeler,Evens}), the detail of the computations is always done for the particular case where only the tilt angle $\theta$ is considered. Since we aim to consider the more general case where tilt and swing angles are used, we will describe the various objects (lines, planes) and the associated computations by using linear algebra tools. 

\subsubsection{The Scheimpflug and the Hinge rules}

In order to explain the Scheimpflug rule, we will refer to the diagram depicted in Figure \ref{scheimpflug}. The  Plane of sharp focus (abbreviated $\mathrm{SFP}$) is determined by a normal vector $\mathbf{n^{SF}}$ and a point $Y$. The position of the optical center $L$ and a vector $\mathbf{n^S}$ normal to the sensor plane (abbreviated $\mathrm{SP}$) are known. The unknowns are the position of the sensor center $S$ and a vector $\mathbf{n^L}$ normal to the lens plane (abbreviated $\mathrm{LP}$). 

The Scheimplug rule stipulates that if $\mathrm{SFP}$ is into focus, then $\mathrm{SP}$, $\mathrm{LP}$ and $\mathrm{SFP}$ necessarily intersect on a common line called the "Scheimpflug Line" (abbreviated $\mathrm{SL}$). The diagram of Figure \ref{scheimpflug}a should help the reader to see that this rule is not sufficient to uniquely determine $\mathbf{n^L}$ and $\mathrm{SP}$, as this plane can be translated toward $\mathbf{n^S}$ if $\mathbf{n^L}$ is changed accordingly. 

The missing constraints are provided by the Hinge rule, which is illustrated in Figure \ref{hinge}. This rule considers two complimentary planes: the front focal plane (abbreviated $\mathrm{FFP}$), which is parallel to $\mathrm{LP}$ and passes through the focal point $F$, and the parallel to sensor lens plane (abbreviated $\mathrm{PSLP}$), which is parallel to $\mathrm{SP}$ and passes through the optical center $L$. The Hinge Rule stipulates that $\mathrm{FFP}$, $\mathrm{PSLP}$ and $\mathrm{SFP}$ must intersect along a common line called the Hinge Line (abbreviated $\mathrm{HL}$). Since $\mathrm{HL}$ is uniquely determined as the intersection of $\mathrm{SFP}$ and $\mathrm{PSLP}$, this allows one to determine $\mathbf{n^L}$, or equivalently the tilt and swing angles, such that $\mathrm{FFP}$ passes through $\mathrm{HL}$ and $F$. Then $\mathrm{SL}$ is uniquely defined as the intersection of $\mathrm{LP}$ and $\mathrm{SFP}$ by the Scheimpflug rule (note that $\mathrm{SL}$ and $\mathrm{HL}$ are parallel by construction). Since $\mathbf{n^S}$ is already known, any point belonging to $\mathrm{SL}$ is sufficient to uniquely define $\mathrm{SP}$. Hence, the determination of tilt and swing angles and position of the rear standard can be summarized as follows:
\begin{enumerate}
\item determination of $\mathrm{HL}$, intersection of $\mathrm{FFP}$ and $\mathrm{SFP}$,
\item determination of tilt and swing angles such that $\mathrm{HL}$ belongs to $\mathrm{FFP}$,
\item determination of $\mathrm{SL}$, intersection of $\mathrm{LP}$ and $\mathrm{SFP}$,
\item translation of $S$ such that $\mathrm{SL}$ belongs to $\mathrm{SP}$.
\end{enumerate}
\subsubsection{Algebraic details of the computations}
\label{computations1}
In this section the origin of the coordinate system is the optical center $L$ and the inner product of two vectors $\mathbf{X}$ and $\mathbf{Y}$ is expressed by using the matrix notation $\scal{\mathbf{X}}{\mathbf{Y}}$. All planes are defined by a unit normal vector and a point in the plane as follows:
\begin{align*}
\mathrm{SP}&=\left\{\mathbf{X}\in\mathbb{R}^3,\;\scal{(\mathbf{X-S})}{\mathbf{n^S}}=0\right\}\\
\mathrm{PSLP}&=\left\{\mathbf{X}\in\mathbb{R}^3,\;\scal{\mathbf{X}}{\mathbf{n^S}}=0\right\},\\
\mathrm{LP}&=\left\{\mathbf{X}\in\mathbb{R}^3,\;\scal{\mathbf{X}}{\mathbf{n^L}}=0\right\},\\
\mathrm{FFP}&=\left\{\mathbf{X}\in\mathbb{R}^3,\;\scal{\mathbf{X}}{\mathbf{n^L}}-f=0\right\},\\
\mathrm{SFP}&=\left\{\mathbf{X}\in\mathbb{R}^3,\;\scal{(\mathbf{X}-\mathbf{Y})}{\mathbf{n^{SF}}}=0\right\},
\end{align*}
where the equation of $\mathrm{FFP}$ takes this particular form because the distance between $L$ and $F$ is equal to the focal length $f$ and we have imposed that $n^L_3>0$. The computations are detailed in the following algorithm:
\begin{algorithm}\rm\mbox{}\label{algo1}%
\begin{itemize}
\item[]{\sl Step 1} : compute the Hinge Line by considering its parametric equation
\begin{align*}
\mathrm{HL}&=\left\{\mathbf{X}\in\mathbb{R}^3,\;\exists\,t\in\mathbb{R},\;\mathbf{X}=\mathbf{W}+t\mathbf{V}\right\},
\end{align*}
where $\mathbf{V}$ is a direction vector and $\mathbf{W}$ is the coordinate vector of an arbitrary point of $\mathrm{HL}$. Since this line is the intersection of $\mathrm{PSLP}$ and $\mathrm{SFP}$, $\mathbf{V}$ is orthogonal to $\mathbf{n^L}$ and $\mathbf{n^{SF}}$. Hence, we can take $\mathbf{V}$ as the cross product 
$$\mathbf{V}=\mathbf{n^{SF}}\times \mathbf{n^S}$$
 and $\mathbf{W}$ as a particular solution (e.g. the solution of minimum norm) of the overdetermined system of equations
\begin{align*}
\scal{\mathbf{W}}{\mathbf{n^S}}&=0,\\
\scal{\mathbf{W}}{\mathbf{n^{SF}}}&=\scal{\mathbf{Y}}{\mathbf{n^{SF}}}.
\end{align*}
\item[] {\sl Step 2} : since $\mathrm{HL}$ belongs to $\mathrm{FFP}$ we have
$$
\scal{(\mathbf{W}+t\mathbf{V})}{\mathbf{n^L}}-f = 0,\quad \forall~t\in\mathbb{R},
$$
hence $\mathbf{n^L}$ verifies the overdetermined system of equations
\begin{align}
\scal{\mathbf{W}}{\mathbf{n^L}}&=f,\label{eq:syshinge1}
\\
\scal{\mathbf{V}}{\mathbf{n^L}}&=0.
\label{eq:syshinge2}
\end{align}
with the constraint $\Vert \mathbf{n^L}\Vert^2=1$. The computation of $\mathbf{n^L}$ can be done by the following two steps:

\begin{enumerate}
\item compute $\widetilde{\mathbf{W}}=\mathbf{V}\times \mathbf{W}$ and $\widetilde{\mathbf{V}}$ the minimum norm solution of system (\ref{eq:syshinge1})-(\ref{eq:syshinge2}), which gives a parametrization $$\mathbf{n^L}=\widetilde{\mathbf{V}}+t\widetilde{\mathbf{W}},$$
of all its solutions, where $t$ is an arbitrary real.
\item determination of $t$ such that $\Vert \mathbf{n^L}\Vert^2=1$: this is done by taking the solution $t$ of the second degree equation
$$
\scal{\widetilde{\mathbf{W}}}{\widetilde{\mathbf{W}}} t^2+2\scal{\widetilde{\mathbf{W}}}{\widetilde{\mathbf{V}}} t+\scal{\widetilde{\mathbf{V}}}{\widetilde{\mathbf{V}}}-1=0,
$$
such that $n^L_3>0$. The tilt and swing angles are then obtained as
$$
\theta_L=-\arcsin n_2^L,\quad \phi_L=-\arcsin \frac{n_1^L}{\cos \theta_L}.
$$
\end{enumerate}
\item[]{\sl Step 3} : since $\mathrm{SL}$ is the intersection of $\mathrm{LP}$ and $\mathrm{SFP}$, the coordinate vector $\mathbf{U}$ of a particular point $U$ on $\mathrm{SL}$ is obtained as the minimum norm solution of the system
\begin{align*}
\scal{\mathbf{U}}{\mathbf{n^L}}&=0,\\
\scal{\mathbf{U}}{\mathbf{n^{SF}}}&=\scal{\mathbf{W}}{\mathbf{n^{SF}}},
\end{align*}
where we have used the fact that $W\in\mathrm{SFP}$.
\item[]{\sl Step 4} : the translation of $S$ can be computed such that $U$ belongs to $\mathrm{SP}$, i.e.
$$
\scal{(\mathbf{U}-\mathbf{S})}{\mathbf{n^S}}=0.
$$
If we only act on the third coordinate of $S$ and leave the two others unchanged, then $S_3$ can be computed as
$$
S_3=\frac{\scal{\mathbf{U}}{\mathbf{n^S}}-S_1n^S_1-S_2n^S_2}{n^S_3}.
$$
\end{itemize}
\end{algorithm}

\begin{remark}\rm
\label{iter}
When we consider a true camera lens, the nodal points $H,H'$ and the front standard center $L$ do not coincide. Hence, the tilt and swing rotations of the front standard modify the actual position of the $\mathrm{PSLP}$ plane. In this case, we use the following iterative fixed point scheme: 

\begin{enumerate}
\item The angles $\phi_L$ and $\theta_L$ are initialized with starting values $\phi_L^0$ and $\theta_L^0$.
\item At iteration $k$, 
\begin{enumerate}
\item \label{a} the position of $\mathrm{PSLP}$ is computed considering $\phi_L^k$ and $\theta_L^k$,
\item \label{b} the resulting Hinge Line is computed, then the position of $\mathrm{FFP}$ and the new values $\phi_L^{k+1}$ and $\theta_L^{k+1}$ are computed.
\end{enumerate}
\end{enumerate}
Point 2 is repeated until convergence of  $\phi_L^k$ and $\theta_L^k$ up to a given tolerance. Generally 3 iterations are sufficient to reach the machine precision. 
\end{remark}

\subsection{Characterization of the depth of field region}

\begin{figure}
\begin{center}
\includegraphics[width=\linewidth]{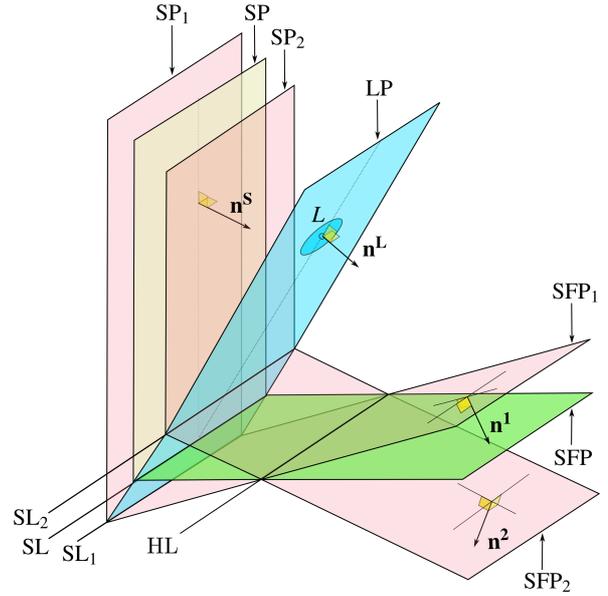}\\
\end{center}
\caption{Position of planes $\mathrm{SFP_1}$ and $\mathrm{SFP_2}$ delimiting the depth of field region.}
\label{scheimpflug-dof}
\end{figure}
As in the previous section, we consider that $L$ and the nodal points $H$ and $H'$ coincide. Moreover, $L$ will be the origin of the global coordinates system.

We consider the configuration depicted in Figure \ref{scheimpflug-dof} where the sharp focus plane $\mathrm{SFP}$, the lens plane $\mathrm{LP}$ and the sensor plane $\mathrm{SP}$ are tied by the Scheimpflug and the Hinge rule. The depth of field can be defined as follows:
\begin{definition}
Let $X$ be a 3D point and $A$ its image through the lens. Let $\mathcal{C}$ be the disk in $\mathrm{LP}$ of center $L$ and diameter $f/N$, where $N$ is called the \textit{f-number}. Let $K$ be the cone of base $\mathcal{C}$ and vertex $A$. The point $X$ is said to lie in the depth of field region if the diameter of the intersection of $\mathrm{SP}$ and $K$ is lower that $c$, the diameter of the so-called \textit{circle of confusion}.
\end{definition}
The common values of $c$, which depend on the magnification from the sensor image to the final image and on its viewing conditions, lie typically between 0.2 mm and 0.01 mm. In the following the value of $c$ is not a degree of freedom but a given input.

If the ellipticity of extended images is neglected, the depth of field region can be shown to be equal to the unbounded wedge delimited by $\mathrm{SFP_1}$ and $\mathrm{SFP_2}$ intersecting at $\mathrm{HL}$, where the corresponding sensor planes $\mathrm{SP_1}$ and $\mathrm{SP_2}$ are tied to $\mathrm{SFP_1}$ and $\mathrm{SFP_2}$ by the Scheimpflug rule. By mentally rotating $\mathrm{SFP}$ around $\mathrm{HL}$, it is easy to see that $\mathrm{SP}$ is translated through $\mathbf{n^S}$ and spans the region between $\mathrm{SP_1}$ and $\mathrm{SP_2}$. The position of $\mathrm{SP_1}$ and $\mathrm{SP_2}$, the f-number $N$ and the diameter of the circle of confusion $c$ are related by the formula
\begin{equation}
\frac{Nc}{f}=\frac{p_1-p_2}{p_1+p_2},
\label{dof}
\end{equation}
where $p_1$, respectively  $p_2$, are the distances between the optical center $L$ and $\mathrm{SP_1}$, respectively $\mathrm{SP_2}$, both measured orthogonally to the optical plane. The distance $p$ between $\mathrm{SP}$ and L can be shown to be equal to
\begin{equation}
p=\frac{2p_1p_2}{p_1+p_2},
\label{dofharm}
\end{equation}
the harmonic mean of $p_1$ and $p_2$. This approximate definition of the depth of field region has been proposed by various authors (see \cite{Wheeler,Bigler}) but when the ellipticity of images is taken into account a complete study can be found in \cite{Evens}. For sake of completeness, we give the justification of formulas (\ref{dof}) and (\ref{dofharm}) in Appendix \ref{dofdemo}. In most practical situations the necessary angle between $\mathrm{SP}$ and $\mathrm{LP}$ is small (less that 10 degrees), so that this approximation is correct.

\begin{remark}
The analysis in Appendix \ref{dofdemo} shows that the ratio $\frac{Nc}{f}$ in equation (\ref{dof}) does not depend on the direction used for measuring the distance between $\mathrm{SP}$, $\mathrm{SP_1}$, $\mathrm{SP_2}$ and $L$. The only condition, in order to take into account the case where $\mathrm{SP}$ and $\mathrm{LP}$ are parallel, is that this direction is not orthogonal to $\mathbf{n^S}$. Hence, by taking the direction given by $\mathbf{n^S}$, we can obtain an equivalent formula to (\ref{dof}). To compute the distances, we need the coordinate vector of two points $U_1$ and $U_2$ on $\mathrm{SP_1}$ and $\mathrm{SP_2}$ respectively. To this purpose we consider Step 3 of Algorithm \ref{algo1} in section \ref{section4}: if $\mathbf{W}$ is the coordinate vector of any point $W$ of $\mathrm{HL}$, each vector $\mathbf{U}^i$ can be obtained as a particular solution of the system 
\begin{align*}
\scal{\mathbf{U}^i}{\mathbf{n^L}}&=0,\\
\scal{\mathbf{U}^i}{\mathbf{n}^i}&=\scal{\mathbf{W}}{\mathbf{n}^i}.
\end{align*}
Since $\mathrm{SP_i}$ can be defined as
$$
\mathrm{SP_i}=\left\{\mathbf{X}\in\mathbb{R}^3,\;\scal{(\mathbf{X}-\mathbf{U}^i)}{\mathbf{n^S}}=0\right\},
$$
and $\Vert\mathbf{n^S}\Vert=1$, the signed distance from $L$ to $\mathrm{SP_i}$ is equal to
$$
d(L,\mathrm{SP_i})=\scal{\mathbf{U_i}}{\mathbf{n^S}}.
$$
So the equivalent formula giving the ratio $\frac{Nc}{f}$ is given by
\begin{align}
\frac{Nc}{f}&=\left\vert\frac{d(L,\mathrm{SP_1})-d(L,\mathrm{SP_2})}{d(L,\mathrm{SP_1})+d(L,\mathrm{SP_2})}\right\vert,\notag\\
&=\left\vert\frac{\scal{(\mathbf{U_1}-\mathbf{U_2})}{\mathbf{n^S}}}{\scal{(\mathbf{U_1}+\mathbf{U_2})}{\mathbf{n^S}}}\right\vert.\label{dof2}
\end{align}
\end{remark}

The above considerations show that for a given orientation of the rear standard given by $\mathbf{n^S}$, if the depth of field wedge is given, then the needed f-number, the tilt and swing angles of the front standard and the translation of the sensor plane, can be determined. The related question that will be addressed in the following is the question: given a set of points $\mathcal{X}=\{{X}^1,\dots,{X}^n\}$, how can we minimize the f-number such that all points of $\mathcal{X}$ lie in the depth of field region?

\subsection{Depth of field optimization with respect to tilt angle}
\label{section5}
We first study the depth of field optimization in two dimensions, because in this particular case all computations can be carried explicitly and a closed form expression is obtained, giving the f-number as a function of front standard tilt angle and of the slope of limiting planes.  First, notice that $N$ has a natural upper bound, since (\ref{dof}) implies that
$$
N \leq \frac{f}{c}.
$$

\subsubsection{Computation of f-number with respect to tilt angle and limiting planes}

\begin{figure}
\begin{center}
\includegraphics[width=\linewidth]{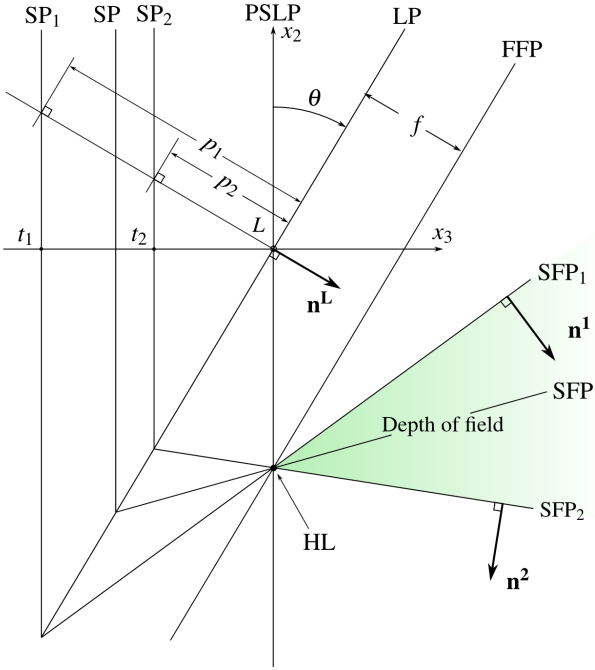}\\
\end{center}
\caption{The depth of field region when only a tilt angle $\theta$ is used.}
\label{dofOptim2D}
\end{figure}
Without loss of generality, we consider that the sensor plane has the normal $\mathbf{n^S}=(0,0,1)^\top$. Let us denote by $\theta$ the tilt angle of the front standard and consider that the swing angle $\phi$ is zero. The lens plane is given by
$$
\mathrm{LP}=\left\{\mathbf{X}\in\mathbb{R}^3,\;\scal{\mathbf{X}}{\mathbf{n^L}}=0\right\},
$$
where
$$
\mathbf{n^L}=(0,-\sin \theta,\cos\theta)^\top,
$$
and any collinear vector to $\mathbf{n^L}\times\mathbf{n^S}$ is a direction vector of the Hinge Line. Hence we can take, independently of $\theta$,  
$$
\mathbf{V}=(1,0,0)^\top
$$
and a parametric equation of $\mathrm{HL}$ is thus given by 
\begin{align*}
\mathrm{HL}&=\left\{\mathbf{X}\in\mathbb{R}^3,\;\exists\,t\in\mathbb{R},\;\mathbf{X}=\mathbf{W}(\theta)+t\mathbf{V}\right\},
\end{align*}
where $\mathbf{W}(\theta)$ is the coordinate vector of a particular point ${W}(\theta)$ on $\mathrm{HL}$, obtained as the minimum norm solution of 
\begin{align*}
\scal{\mathbf{W}(\theta)}{\mathbf{n^S}}&=0,\\
\scal{\mathbf{W}(\theta)}{\mathbf{n^L}}&=f.
\end{align*}
Straightforward computations show that for $\theta\neq 0$
$$\mathbf{W}(\theta)=\left(\begin{array}{c}0\\-\frac{f}{\sin\theta}\\0\end{array}\right).$$

Consider, as depicted in Figure \ref{dofOptim2D}, the two sharp focus planes $\mathrm{SFP_1}$ and $\mathrm{SFP_2}$ passing through $\mathrm{HL}$, with normals $\mathbf{n^1}=(0,-1,a_1)^\top$ and $\mathbf{n^2}=(0,-1,a_2)^\top$,
\begin{align*}
\mathrm{SFP_1}&=\left\{\mathbf{X}\in\mathbb{R}^3,\;\scal{(\mathbf{X}-\mathbf{W}(\theta))}{\mathbf{n^1}}=0\right\},\\
\mathrm{SFP_2}&=\left\{\mathbf{X}\in\mathbb{R}^3,\;\scal{(\mathbf{X}-\mathbf{W}(\theta))}{\mathbf{n^2}}=0\right\}.
\end{align*}
The two corresponding sensor planes $\mathrm{SP_i}$ are given by following Steps 3-4 in Algorithm 1 by 
\begin{align*}
\mathrm{SFP_i}=\left\{\mathbf{X}\in\mathbb{R}^3,\;X_3=t_i\right\},\;i=1,2,
\end{align*}
where
$$
t_i=\frac{f}{a_i\sin\theta-\cos\theta}, \;i=1,2.
$$
Using equation (\ref{dof2}) the corresponding f-number is equal to 
\begin{equation}
N(\theta,\mathbf{a})=\left\vert\frac{t_1-t_2}{t_1+t_2}\right\vert\left(\frac{f}{c}\right),
\end{equation}
where we have used the notation $\mathbf{a}=(a_1,a_2)$. Finally we have, for $\theta\neq 0$
\begin{equation}
N(\theta,\mathbf{a})=\operatorname{sign}\theta\frac{(a_1-a_2)\sin\theta}{2\cos \theta-(a_1+a_2)\sin\theta}\left(\frac{f}{c}\right).
\label{formulaN}
\end{equation}
\begin{remark}\rm When $\theta=0$, then $\mathrm{SFP}$ does not intersect $\mathrm{PSLP}$ and the depth of field region is included between two parallel planes $\mathrm{SFP_i}$ given by 
\begin{align*}
\mathrm{SFP_i}=\left\{\mathbf{X}\in\mathbb{R}^3,\;X_3=z_i\right\},\;i=1,2,
\end{align*}
where $z_1$ and $z_2$ depend on the f-number and on the position of $\mathrm{SFP}$. One can show by using the thin lens equation and equation (\ref{dof}) that the corresponding f-number is equal to
\begin{equation}
N_0(z_1,z_2)=\frac{\left\vert \frac{1}{z_2}-\frac{1}{z_1}\right\vert}{\frac{1}{z_1}+\frac{1}{z_2}-\frac{2}{f}}\left(\frac{f}{c}\right).
\end{equation}
\end{remark}

\subsubsection{Theoretical results for the optimization problem}

Without loss of generality, we will consider a set of only three non-aligned points $\mathcal{X}=\{{X}^1,{X}^2,{X}^3\}$, which have to be within the depth of field region with minimal f-number and we denote by $\mathbf{X}^1,\mathbf{X}^2,\mathbf{X}^3$ their respective coordinate vectors.

The corresponding optimization problem can be stated as follows: find
\begin{equation}
(\theta^*,\mathbf{a}^*)=\arg\min_{\begin{array}{c}\theta\in\mathbb{R}\\ a\in \mathcal{A(\theta)}\end{array}} N(\theta,\mathbf{a}),
\end{equation}
where for a given $\theta$ the set $\mathcal{A(\theta)}$ is defined by the inequalities
\begin{align}
\scal{(\mathbf{X}^i-\mathbf{W}(\theta))}{\mathbf{n^1}}&\geq 0,~i=1,2,3,\label{cf1}\\
\scal{(\mathbf{X}^i-\mathbf{W}(\theta))}{\mathbf{n^2}}&\leq 0,~i=1,2,3,\label{cf2}
\end{align}
meaning that ${X}^1,{X}^2,{X}^3$ are respectively under $\mathrm{SFP_1}$ and above $\mathrm{SFP_2}$, and by the inequalities
\begin{align}
-\scal{(\mathbf{X}^i-\mathbf{W}(\theta))}{\mathbf{n^L}}+f&\leq 0,~i=1,2,3,\label{c1}
\end{align}
meaning that ${X}^1,{X}^2,{X}^3$ are in front of $\mathrm{FFP}$.
\begin{remark}
\rm Points behind the focal plane cannot be in focus, and the constraints (\ref{c1}) are just expressing this practical impossibility. However, we have to notice that when one of these constraints is active, we can show that $a_1=\cot\theta$ or $a_2=\cot\theta$ so that $N(\mathbf{a},\theta)$ reaches its upper bound $\frac{f}{c}$. Moreover, we have to eliminate the degenerate case where the points ${X}^i$ are such that there exists two active constraints in (\ref{c1}): in this case, there exists a unique admissible pair $(\mathbf{a},\theta)$ and the problem has no practical interest. To this purpose, we can suppose that $X^i_3>f$ for $i=1,2,3$.
\end{remark}

For $\theta\neq 0$ the gradient of $N(\mathbf{a},\theta)$ is equal to
$$
\mathbf{\nabla}(\mathbf{a},\theta)=\frac{2\operatorname{sign}\theta}{(2\cot \theta -(a_1+a_2))^2}\left(\begin{array}{c}\cot\theta-a_2\\-\cot\theta+a_1\\\frac{a_1-a_2}{\sin^2\theta}\end{array}\right)
$$
and cannot vanish since $a_1=a_2$ is not possible because it would mean that ${X}^1,{X}^2,{X}^3$ are aligned. This implies that $\mathbf{a}$ lies on the boundary of $\mathcal{A}(\theta)$ and we have the following intuitive result (the proof is given in Appendix \ref{proof:prop1}):
\begin{proposition}\label{prop1}Suppose that  $X^i_3>f$, for $i=1,2,3$. Then when $N(\mathbf{a},\theta)$ reaches its minimum, there exists $i_1$, $i_2$ with $i_1\neq i_2$ such that
\begin{align*}
\scal{(\mathbf{X}^{i_1}-\mathbf{W}(\theta))}{\mathbf{n^1}}&= 0,\\
\scal{(\mathbf{X}^{i_2}-\mathbf{W}(\theta))}{\mathbf{n^2}}&= 0.
\end{align*}
\end{proposition}
\begin{remark}
\rm The above result shows that at least two points touch the depth of field limiting planes $\mathrm{SFP_1}$ and $\mathrm{SFP_2}$ when the f-number is minimal. In the following, we will show that the three points ${X}^1,{X}^2$ and ${X}^3$ are necessarily in contact with one of the limiting planes (the proof is given in Appendix \ref{proof:prop2}):
\end{remark}
\begin{proposition}\label{prop2}Suppose that the vertices $\{{X}^i\}_{i=1\dots 3}$ verify the condition
\begin{equation}
\frac{\Vert \mathbf{X}^{i}\times \mathbf{X}^{j}\Vert}{\Vert \mathbf{X}^{i}-\mathbf{X}^{j} \Vert}>f,~~i\neq j.
\label{conddet}
\end{equation}
Then $N(\mathbf{a},\theta)$ reaches its minimum when all vertices are in contact with the limiting planes.
\end{proposition}
\begin{remark}\rm
If $\theta$ is small, then $N(\theta,\mathbf{a})$ in (\ref{formulaN}) can be approximated by 
\begin{equation}
\tilde N(\theta,\mathbf{a})=\operatorname{sign}\theta{(a_1-a_2)\sin\theta}\left(\frac{f}{2c}\right),
\label{approxN}
\end{equation}
and the proof of Proposition \ref{prop2} is considerably simplified: the same result holds with the weaker condition
\begin{equation}
\Vert \mathbf{X}^{i}\times \mathbf{X}^{j}\Vert> 0.
\end{equation}
In fact, an approximate way of specifying the depth of field region using the \textit{hyperfocal distance}, proposed in \cite{Merklinger}, leads to the same approximation of $N(\theta,\mathbf{a})$, under the {\em a priori} hypothesis of small $\theta$ and distant objects, i.e. $X^i_3\gg f$. This remark is clarified in Appendix \ref{appendix1}.
\end{remark}
We will illustrate the theoretical result by considering sets of 3 points. For a set with more than 3 vertices (but being equal to the vertices of the convex hull of $\mathcal{X}$), the determination of the optimal solution is purely combinatorial, since it is enough to enumerate all admissible situations where two points are in contact with one plane, and a third one with the other. The value $\theta=0$ also has to be considered because it can be a critical value if the object has a vertical edge. We will also give an Example which violates condition (\ref{conddet}) and where $N(\mathbf{a},\theta)$ reaches its minimum when only two vertices are in contact with the limiting planes.
\subsubsection{Numerical results}
\label{2dnum}
In this section, we will consider the following function, defined for $\theta\neq 0$
$$
n(\theta)=\min_{a\in\mathcal{A}(\theta)}N(\mathbf{a},\theta).
$$
Finding the minimum of this function allows one to solve the original constrained optimization problem, but considering the results of the previous section, $n(\theta)$ is non-differentiable. In fact, the values of $\theta$ for which $n(\theta)$ is not differentiable correspond to the situations where 3 points are in contact with the limiting planes. We extend $n(\theta)$ by continuity for $\theta=0$ by defining
$$
n(0)=\frac{\frac{1}{z_2}-\frac{1}{z_1}}{\frac{1}{z_1}+\frac{1}{z_2}-\frac{2}{f}}\left(\frac{f}{c}\right),
$$
where 
$$
z_1=\max_{i=1,2,3}X^i_3,~~z_2=\min_{i=1,2,3}X^i_3.
$$
This formula can be directly obtained by using conjugation formulas or by taking $\theta=0$ in equation (\ref{dofproof}).

\begin{example}\rm\mbox{}
\begin{figure}
\begin{center}
\includegraphics[width=\linewidth]{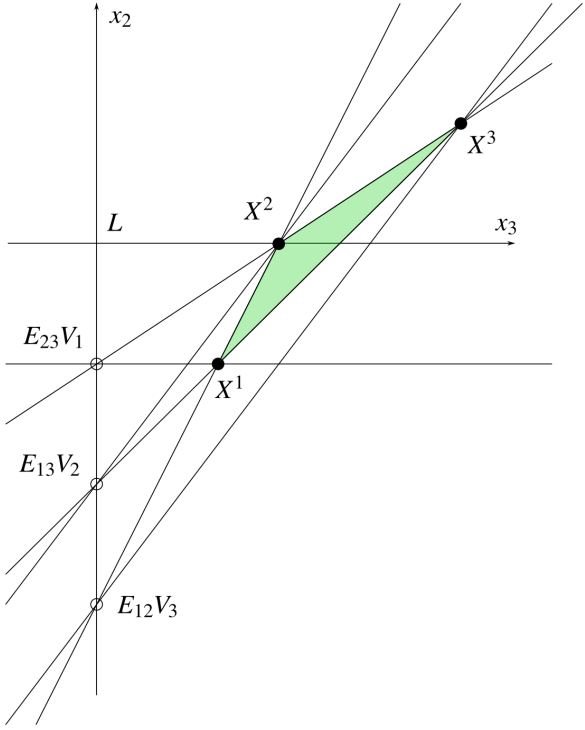}
\end{center}
\caption{Enumeration of the 3 candidates configurations for Example 1.}
\label{figexample1}
\end{figure}
\begin{figure}
\begin{center}
\begin{tabular}{c}
\includegraphics[width=0.95\linewidth]{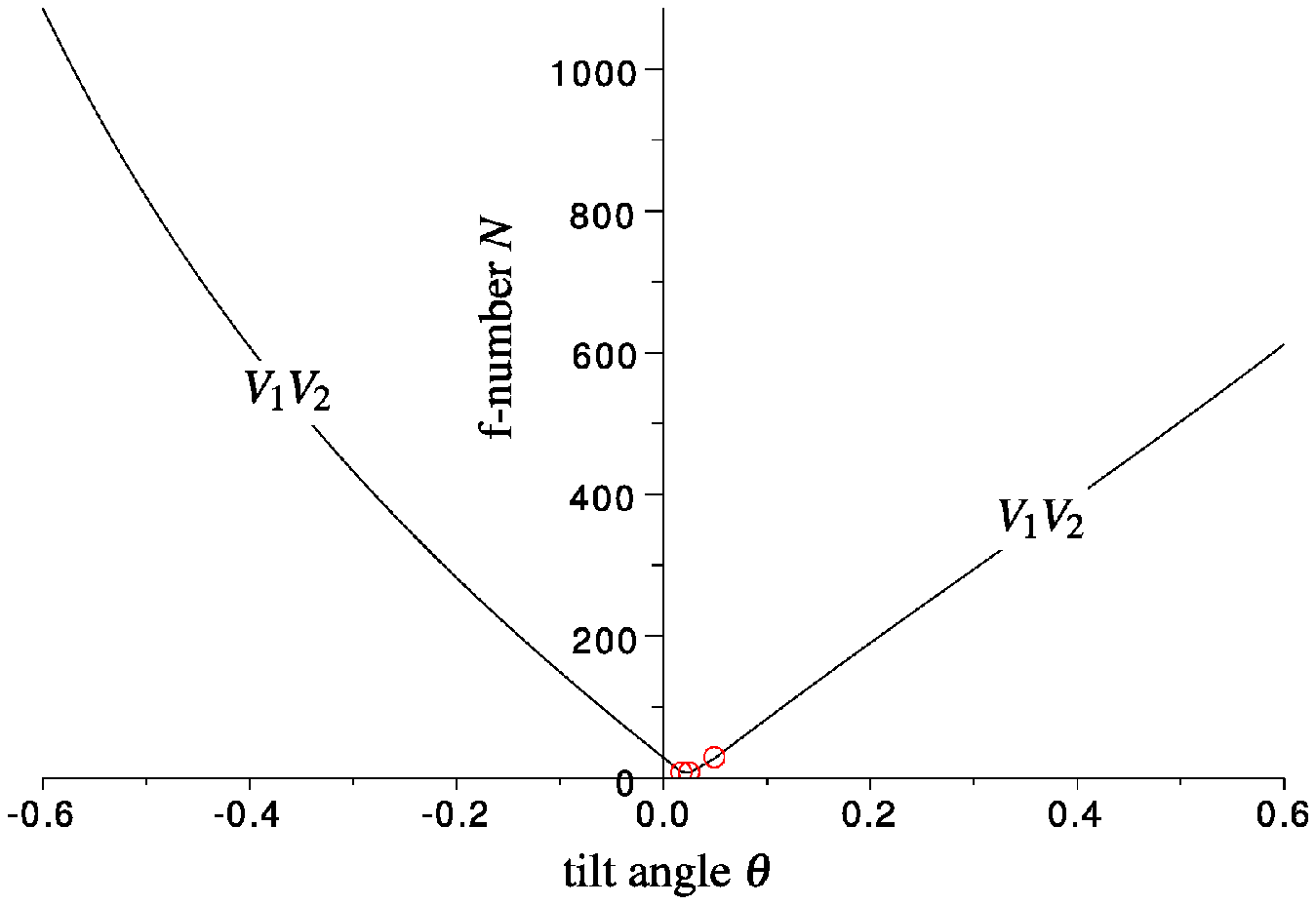}\\(a)\\
\includegraphics[width=0.95\linewidth]{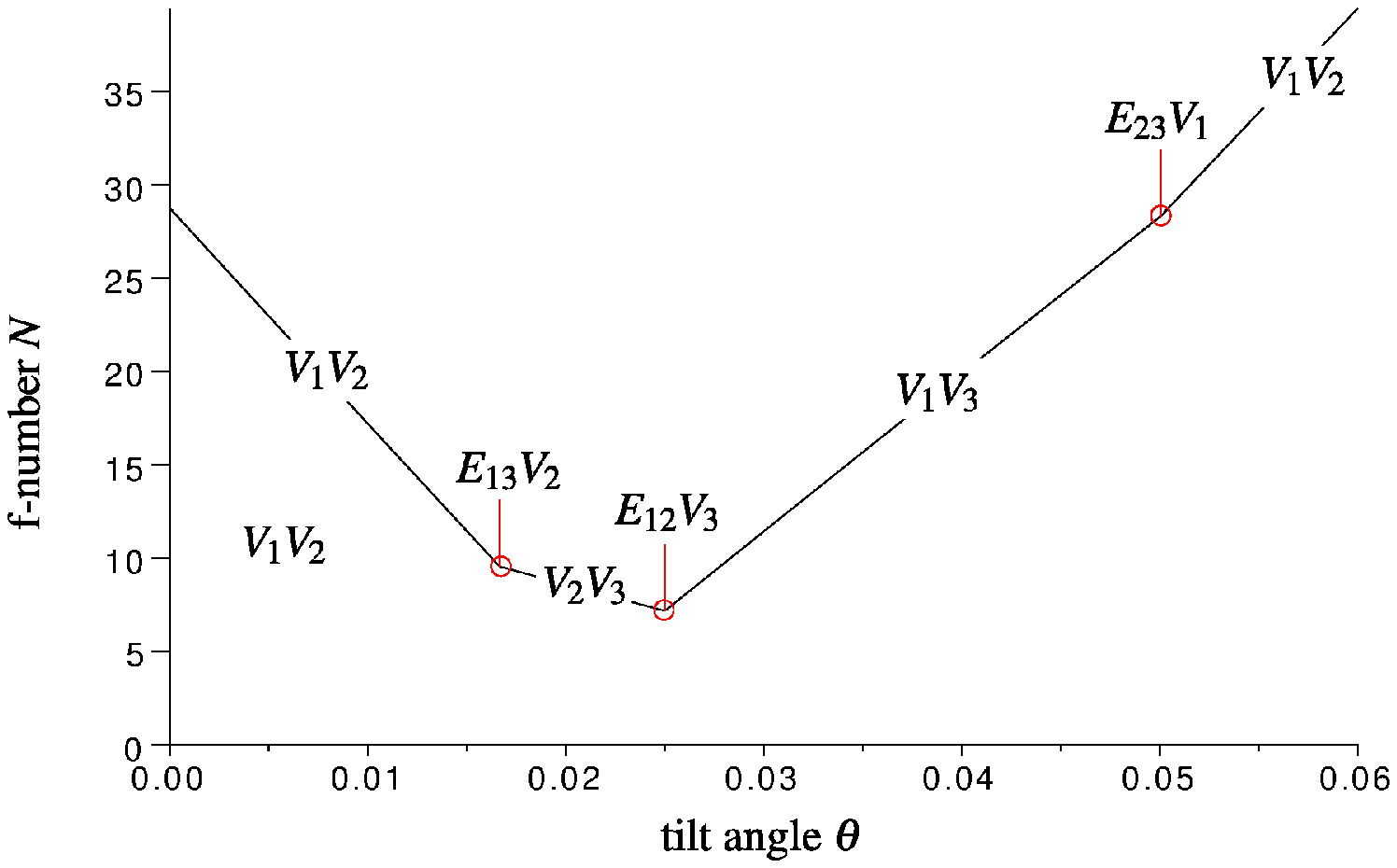}\\(b)\\
\end{tabular}
\end{center}
   \caption{Example 1 graph of $n(\theta)$ for $\theta\in[-0.6,0.6]$ (a) and  $\theta\in[0,0.06]$ (b). Labels give the different types of contact.}
\label{figexample1curves}
\end{figure}
We consider a lens with focal length $f=5.10^{-2}\mathrm{m}$ and a confusion circle $c=3.10^{-5}m$ (commonly used value for 24x36 cameras). The vertices have coordinates 
$$
\mathbf{X}^1=\left(\begin{array}{r}0\\-1\\1\end{array}\right),~\mathbf{X}^2=\left(\begin{array}{r}0\\3\\1\end{array}\right),~ \mathbf{X}^3=\left(\begin{array}{r}0\\0\\1.5\end{array}\right).
$$
Figure \ref{figexample1} shows the desired depth of field region (dashed zone) and the three candidates hinge lines corresponding to contacts $E_{ij}V_{k}$:
\begin{itemize} 
\item $E_{12}V_3$, obtained when the edge $[{X}^1,{X}^2]$ is in contact with $\mathrm{SFP_1}$ and  vertex ${X}^3$ with $\mathrm{SFP_2}$.
\item $E_{13}V_2$, obtained when the edge $[{X}^1,{X}^3]$ is in contact with $\mathrm{SFP_2}$ and vertex ${X}^2$ with $\mathrm{SFP_1}$.
\item $E_{23}V_1$, obtained when the edge $[{X}^2,{X}^3]$ is in contact with $\mathrm{SFP_1}$ and vertex ${X}^1$ with $\mathrm{SFP_2}$.
\end{itemize}
The associated values of $\theta$ and $n(\theta)$ are given in Table \ref{table2D},
\begin{table}
\begin{center}
\begin{tabular}{|c|c|c|c|}
\hline Contact & $\theta$ &  $n(\theta)$\\
\hline $E_{12}V_{3}$ & 0.0166674 & 9.52\\
\hline $E_{13}V_{2}$ & 0.025002 & 7.16\\
\hline $E_{23}V_{1}$ & 0.0500209 & 28.28\\
\hline
\end{tabular}
\end{center}
\caption{Value of $\theta$ for each possible optimal contact and corresponding f-number $n(\theta)$ for Example 1.}
\label{table2D}
\end{table}
which shows that contact $E_{13}V_{2}$ seems to give the minimum f-number. Condition (\ref{conddet}) is verified, since
$$
\frac{\Vert \mathbf{X}^1\times \mathbf{X}^2\Vert}{\Vert \mathbf{X}^1-\mathbf{X}^2\Vert}=1,~\frac{\Vert \mathbf{X}^1\times \mathbf{X}^3\Vert}{\Vert \mathbf{X}^1-\mathbf{X}^3\Vert}=1.34,$$
$$
\frac{\Vert \mathbf{X}^2\times \mathbf{X}^3\Vert}{\Vert \mathbf{X}^2-\mathbf{X}^3\Vert}=1.48.
$$
Hence, the derivative of $n(\theta)$ cannot vanish and the minimum f-number is necessarily reached for $\theta^*=0.025002$.

We can confirm this by considering the graph of $n(\theta)$ depicted on Figure \ref{figexample1curves}. In the zone of interest, the function $n(\theta)$ is almost piecewise affine. For each point in the interior of curved segments of the graph, the value of $\theta$ is such that the contact of $\mathcal{X}$ with the limiting planes is of type $V_iV_j$. Clearly, the derivative of $n(\theta)$ does not vanish. The possible optimal values of $\theta$, corresponding to contacts of type $E_ijV_k$, are the abscissa of angular points of the graph, marked with red dots. The graph confirms that the minimal value of $n(\theta)$ is reached for contact $E_{12}V_3$.

The minimal f-number is equal to $n(\theta^*)=7.16$. By comparison, the f-number without tilt optimization is $n(0)=28.74$. This example highlights the important gain in terms of f-number reduction with the optimized tilt angle.

\label{example1_2d}
\end{example}

\begin{example}\rm\mbox{}
\begin{table}
\begin{center}
\begin{tabular}{|c|c|c|c|}
\hline Contact & $\theta$ &  $n(\theta)$\\
\hline $E_{12}V_{3}$ & 0.185269 & 29.49\\
\hline $E_{23}V_{1}$ & 0.100419 & 83.67\\
\hline $E_{13}V_{2}$ & 0.235825 & 47.04\\
\hline
\end{tabular}
\end{center}
\caption{Value of $\theta$ for each possible optimal contact and corresponding f-number $n(\theta)$ for Example 2.}
\label{table2D2}
\end{table}
\begin{figure}
\begin{center}
\begin{tabular}{c}
\includegraphics[width=0.95\linewidth]{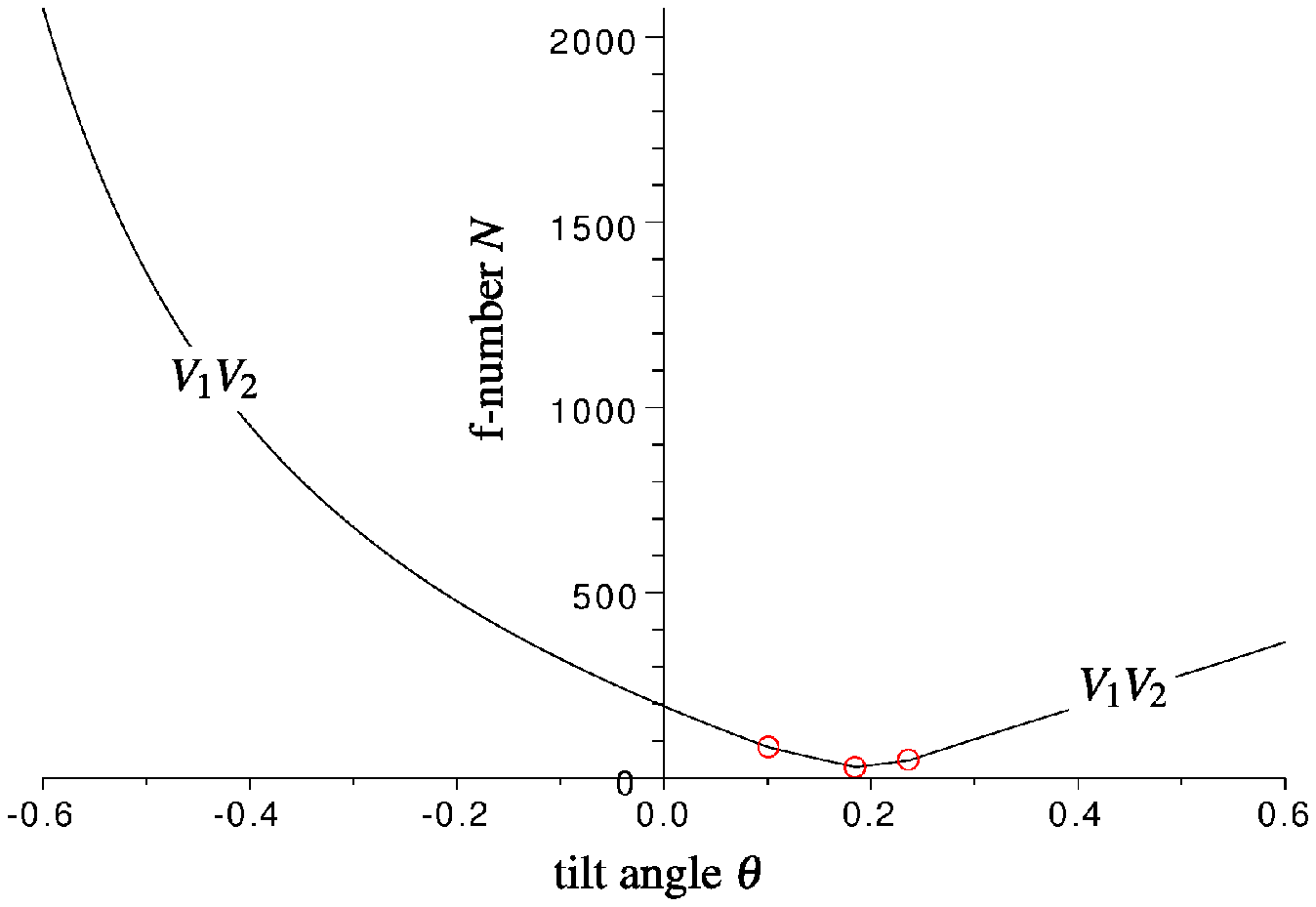}\\(a)
\\\includegraphics[width=0.95\linewidth]{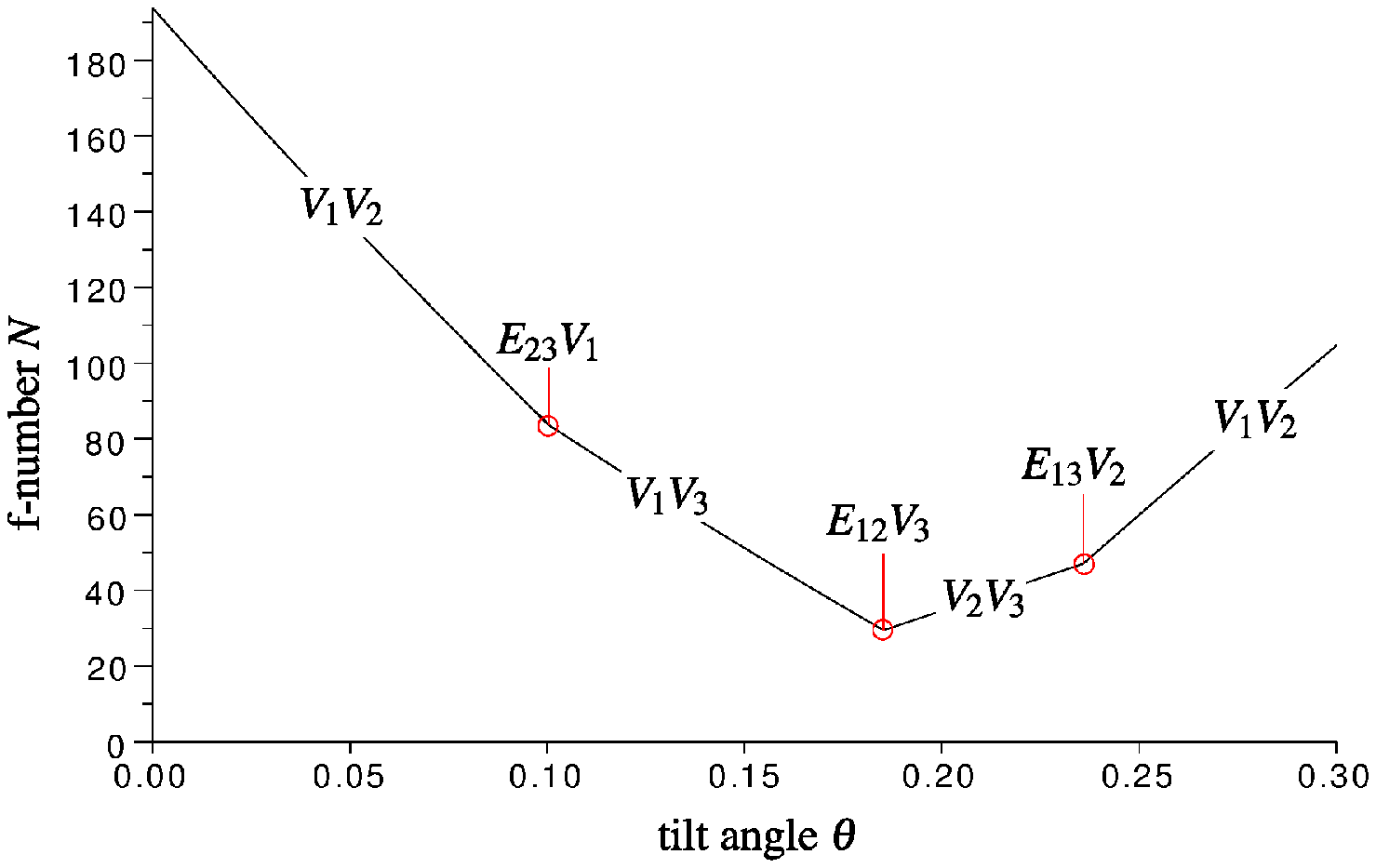}\\(b)
\end{tabular}
\end{center}
   \caption{Example 2 graph of $n(\theta)$ for $\theta\in[-0.6,0.6]$ (a) and  $\theta\in[0,0.3]$ (b). Labels give the different types of contact.}
\label{figexample2curves}
\end{figure}
We consider the same lens and confusion circle as in Example 1 ($f=5.10^{-2}\mathrm{m}$, $c=3.10^{-5}m$) but the vertices have coordinates 
$$
\mathbf{X}^1=\left(\begin{array}{r}0\\-0.1\\0.12\end{array}\right),~\mathbf{X}^2=\left(\begin{array}{r}0\\0\\0.19\end{array}\right),~ \mathbf{X}^3=\left(\begin{array}{r}0\\-0.0525\\0.17\end{array}\right).
$$
The object is almost ten times smaller than the object of the previous example (it has the size of a small pen), but it is also ten times closer: this is a typical close up configuration. The values of $\theta$ and $n(\theta)$ associated to contacts of type $E_{ij}V_k$ are given in Table \ref{table2D2} which shows that contact $E_{13}V_{2}$ seems to give the minimum f-number. We have
$$
\frac{\Vert \mathbf{X}^1\times \mathbf{X}^2\Vert}{\Vert \mathbf{X}^1-\mathbf{X}^2\Vert}=0.16,~\frac{\Vert \mathbf{X}^1\times \mathbf{X}^3\Vert}{\Vert \mathbf{X}^1-\mathbf{X}^3\Vert}=0.16,$$
$$
\frac{\Vert \mathbf{X}^2\times \mathbf{X}^3\Vert}{\Vert \mathbf{X}^2-\mathbf{X}^3\Vert}=0.1775527,
$$
showing that condition (\ref{conddet}) is still verified, even if the values are smaller than the values of Example 1. Hence, the derivative of $n(\theta)$ cannot vanish and the minimum f-number is reached for $\theta^*=0.185269$. 

We can confirm this by considering the graph of $n(\theta)$ depicted on Figure \ref{figexample2curves}. As in Example 1, the derivative of $n(\theta)$ does not vanish and the graph confirms that the minimal value of $n(\theta)$ is reached for contact $E_{12}V_3$.

The minimal f-number is equal to $n(\theta^*)=29.49$. By comparison, the f-number without tilt optimization is $n(0)=193.79$. Such a large value gives an aperture of diameter $0.26$mm, almost equivalent to a pin hole ! Since the maximum f-number of view camera lenses is never larger than 64, the object cannot be in focus without using tilt.

This example also shows that Proposition \ref{prop2} is still valid even if the object is close to the lens and the obtained optimal tilt angle $0.185269$ (10.61 degrees) is large.
\end{example}

\begin{example}\label{example3}\rm\mbox{}
\begin{figure}
\begin{center}
\begin{tabular}{c}
\includegraphics[width=0.95\linewidth]{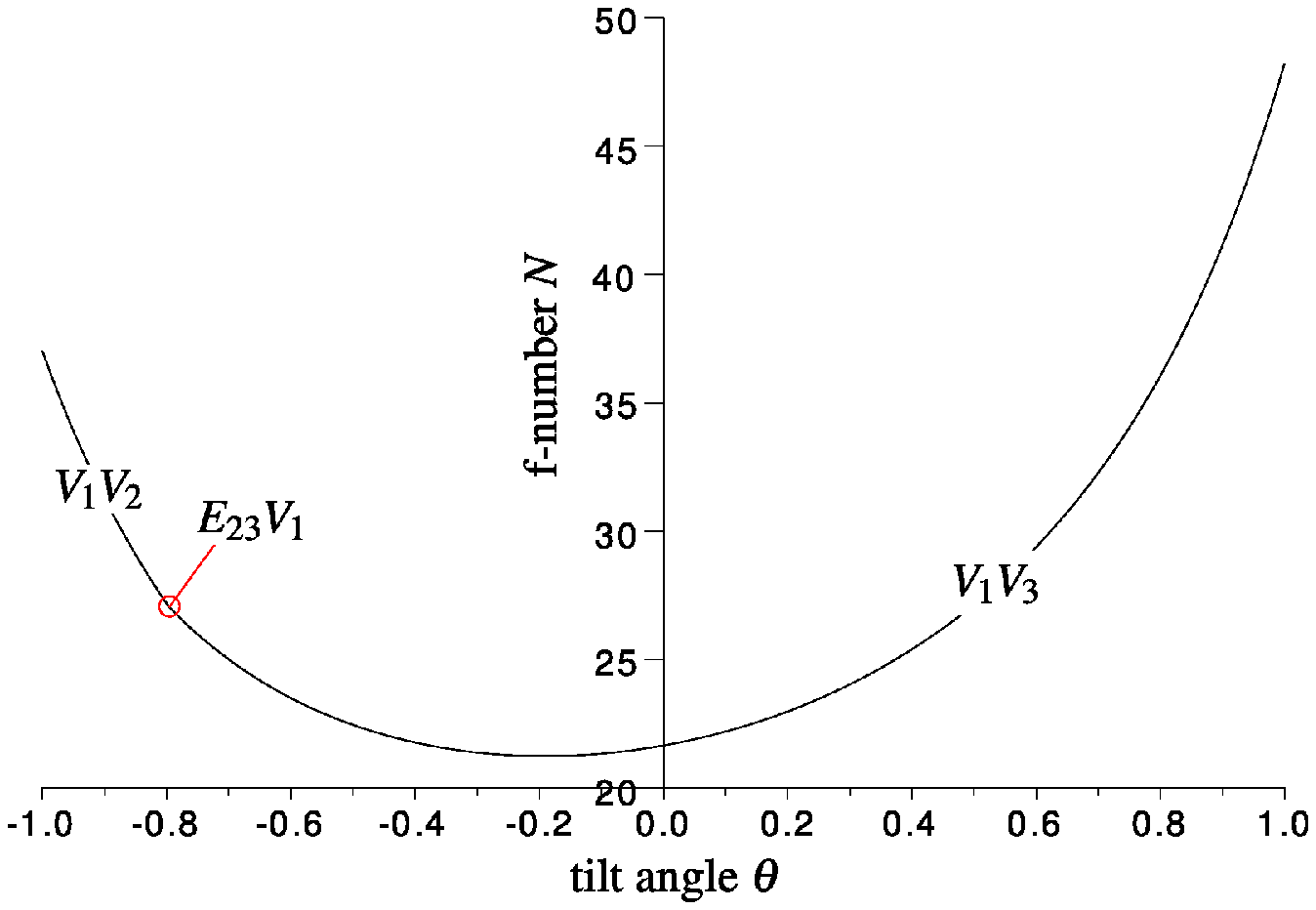}\\(a)\\
\includegraphics[width=0.95\linewidth]{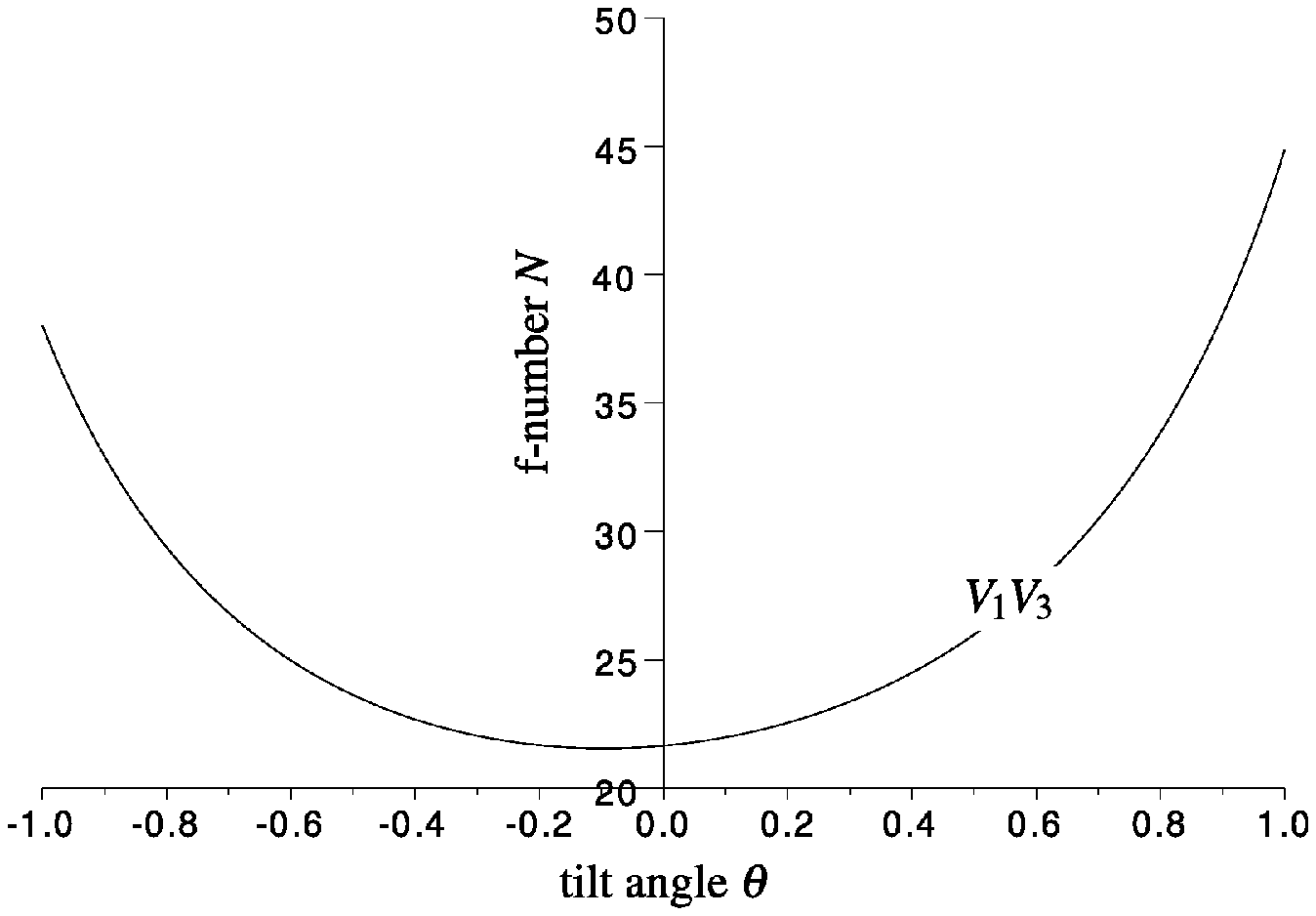}\\(b)
\end{tabular}
\end{center}
   \caption{Example 3 graphs of $n(\theta)$ for $\theta\in[-1,1]$.}
\label{figexample3curves}
\end{figure}
In this example we consider an extremely flat triangle which is almost aligned with the optical center. From the photographer's point of view, this is an unrealistic case. We consider the same optical parameters as before and consider the following vertices
$$
\mathbf{X}^1=\left(\begin{array}{r}0\\0\\1\end{array}\right),~\mathbf{X}^2=\left(\begin{array}{r}0\\h\\1.5\end{array}\right),~ \mathbf{X}^3=\left(\begin{array}{r}0\\-h\\2\end{array}\right),
$$
for $h=0.01$ and we have
$$
\frac{\Vert \mathbf{X}^1\times \mathbf{X}^2\Vert}{\Vert \mathbf{X}^1-\mathbf{X}^2\Vert}=0.02,~\frac{\Vert \mathbf{X}^1\times \mathbf{X}^3\Vert}{\Vert \mathbf{X}^1-\mathbf{X}^3\Vert}=0.01,$$
$$
\frac{\Vert \mathbf{X}^2\times \mathbf{X}^3\Vert}{\Vert \mathbf{X}^2-\mathbf{X}^3\Vert}=0.07,
$$
hence condition (\ref{conddet}) is violated in configurations $V_1V_2$ and $V_1V_3$. Since this condition is sufficient, the minimum value of $n(\theta)$ could still occur for a contact of type $E_{ij}V_k$. However, we can see in Figure \ref{figexample3curves}a that $n(\theta)$ has a minimum at a differentiable point in configuration $V_1V_3$. 
\end{example}
\begin{remark}\rm
Condition (\ref{conddet}) is not necessary: in the proof of Proposition \ref{prop2} (given in Appendix \ref{proof:prop2}), it can be seen that (\ref{conddet}) is a condition ensuring that the polynomial $p(\theta)$ defined by equation (\ref{ptheta}) has no root in $[-\pi/2,\pi/2]$. However, the relevant interval is smaller because  ${X}^{i_1}$ and ${X}^{i_2}$ do not stay in contact with the limiting planes for all values of $\theta$ in $[-\pi/2,\pi/2]$, and because all vertices must be in front of the focal plane. Anyway,  it is always possible to construct absolutely unrealistic configurations. For example, when we consider the above vertices with $h=0.005$, then no edge is in contact with the limiting planes and vertices ${X}^1$ and ${X}^3$ stay in contact with the limiting planes for all admissible values of $\theta$. The corresponding graph of $n(\theta)$ is given in Figure \ref{figexample3curves}b.
\end{remark}

\subsection{Depth of field optimization with respect to tilt and swing angles}

As in the previous section, we consider the case of a thin lens and where the optical and sensor centers coincide respectively with the front and rear standard rotation centers. Without loss of generality, we consider that the sensor plane has the normal $\mathbf{n^S}=(0,0,1)^\top$. The lens plane is given by
$$
\mathrm{LP}=\left\{\mathbf{X}\in\mathbb{R}^3,\;\scal{\mathbf{X}}{\mathbf{n^L}}=0\right\},
$$
where
$$
\mathbf{n^L}=(-\sin\phi\cos\theta,-\sin\theta,\cos\phi\cos\theta)^\top.
$$
A parametric equation of $\mathrm{HL}$ is given by 
\begin{align*}
\mathrm{HL}&=\left\{\mathbf{X}\in\mathbb{R}^3,\;\exists\,t\in\mathbb{R},\;\mathbf{X}=\mathbf{W}(\theta,\phi)+t\mathbf{V}(\theta,\phi)\right\},
\end{align*}
where the direction vector is given by 
$$
\mathbf{V}(\theta,\phi)=\mathbf{n^L}\times\mathbf{n^S}=(-\sin\theta,\sin\phi\cos\theta,0)^\top,
$$
and $\mathbf{W}(\theta,\phi)$ is the coordinate vector of a particular point ${W}(\theta,\phi)$ on $\mathrm{HL}$, obtained as the minimum norm solution of 
\begin{align*}
\scal{\mathbf{W}(\theta)}{\mathbf{n^S}}&=0,\\
\scal{\mathbf{W}(\theta)}{\mathbf{n^L}}&=f.
\end{align*}
Consider, as depicted in Figure \ref{scheimpflug-dof}, the two planes of sharp focus $\mathrm{SFP_1}$ and $\mathrm{SFP_2}$ intersecting at HL, with normals $\mathbf{n^1}$ and $\mathbf{n^1}$ respectively. The point $W(\theta,\phi)$ belongs to $\mathrm{SFP_1}$ and $\mathrm{SFP_2}$ and any direction vector of  $\mathrm{SFP_1}\cap\mathrm{SFP_2}$ is collinear to $\mathbf{V}(\theta,\phi)$. Hence, we have
\begin{align*}
\mathrm{SFP_1}&=\left\{\mathbf{X}\in\mathbb{R}^3,\;\scal{(\mathbf{X}-\mathbf{W}(\theta,\phi))}{\mathbf{n^1}}=0\right\},\\
\mathrm{SFP_2}&=\left\{\mathbf{X}\in\mathbb{R}^3,\;\scal{(\mathbf{X}-\mathbf{W}(\theta,\phi))}{\mathbf{n^2}}=0\right\},
\end{align*}
and 
\begin{equation}
(\mathbf{n^1}\times\mathbf{n^2})\times \mathbf{V}(\theta,\phi)=0.
\label{eqn1n2HL}
\end{equation}
Using equation (\ref{dof2}) the f-number
is equal to
\begin{equation}
N(\theta,\phi,\mathbf{n^1},\mathbf{n^2})=\frac{f}{c}\left\vert\frac{\scal{(\mathbf{U}^1-\mathbf{U}^2)}{\mathbf{n^S}}}{\scal{(\mathbf{U}^1+\mathbf{U}^2)}{\mathbf{n^S}}}\right\vert,\label{fnumber3d}
\end{equation}
where each coordinate vector $\mathbf{U}^i$ of point $U^i$, for $i=1,2$, is obtained as a particular solution of the following system:
\begin{align*}
\scal{\mathbf{U}^i}{\mathbf{n^L}}&=0,\\
\scal{\mathbf{U}^i}{\mathbf{n}^i}&=\scal{\mathbf{W}(\theta,\phi)}{\mathbf{n}^i}.
\end{align*}

\begin{remark}\rm\mbox{}\begin{figure}
\begin{center}
\includegraphics[width=\linewidth]{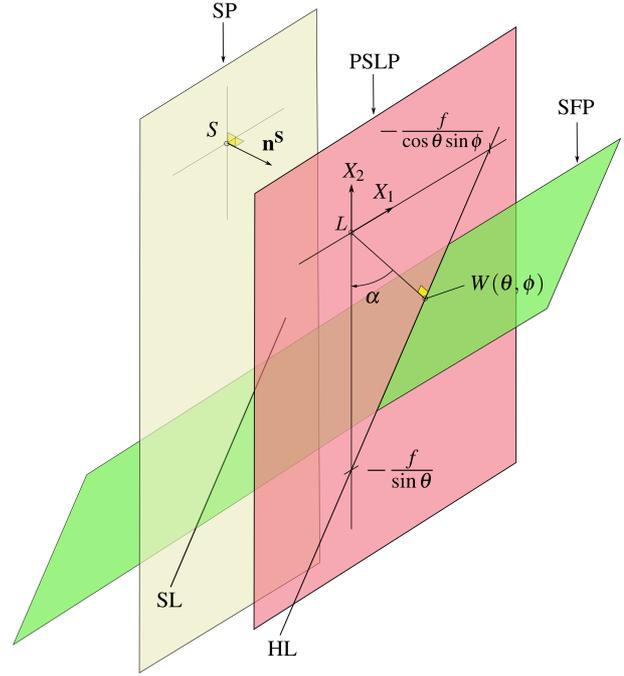}
\end{center}
\caption{Position of the hinge line when both tilt and swing are used and $\mathbf{n^S}=(0,0,1)^\top$. The LP and FFP planes have not been represented for reasons of readability.}
\label{hinge3d}
\end{figure}
When $\mathbf{n^S}=(0,0,1)^\top$ the intersections of $\mathrm{HL}$ with the $(L,X_1)$ and the $(L,X_2)$ axes can be determined, as depicted in Figure \ref{hinge3d}. In this case, it is easy to show that the coordinates of the minimum norm $\mathbf{W}(\theta,\phi)$ are given by
$$
\mathbf{W}(\theta,\phi)=-\frac{f}{\sin^2\phi\cos^2\theta+\sin^2\theta}(\sin\phi\cos\theta,\sin\theta,0)^\top.
$$
We note that up to a rotation of axis $(0,0,1)^\top$ and angle $\alpha$, we recover a configuration where only a tilt angle $\psi$ is used, where these two angles are respectively defined by
\begin{align*}
\sin\psi&=\operatorname{sign}\theta\sqrt{\sin^2\phi\cos^2\theta+\sin^2\theta},\\
\sin\alpha&=\frac{\sin\phi\cos\theta}{\sin\psi}.
\end{align*}
\label{remrot}
\end{remark}
\subsubsection{Optimization problem}
Consider a set 4 non coplanar points $\mathcal{X}=\{{X}^i\}_{i=1\dots 4}$, which have to be within the depth of field region with minimal f-number and denote by $\{\mathbf{X}^i\}_{i=1\dots 4}$ their respective coordinate vectors. The optimization problem can be stated as follows: find
\begin{equation}
(\theta^*,\phi^*,\mathbf{n^1}^*,\mathbf{n^2}^*)=\arg\min_{\theta,\phi,\mathbf{n^1},\mathbf{n^2}} N(\theta,\phi,\mathbf{n^1},\mathbf{n^2}),
\label{optim}
\end{equation}
where the minimum is taken for $\theta,\phi,\mathbf{n^1}$ and $\mathbf{n^2}$ such that equation (\ref{eqn1n2HL}) is verified and such that the points $\{{X}^i\}_{i=1\dots 4}$ lie between $\mathrm{SFP_1}$ and $\mathrm{SFP_2}$. This last set of constraints can be expressed in a way similar to equations (\ref{cf1})-(\ref{cf2}).

\subsubsection{Analysis of configurations}

Consider the function $n(\theta,\phi)$ defined by 
\begin{equation}
n(\theta,\phi)=\min_{\mathbf{n^1},\mathbf{n^2}} N(\theta,\phi,\mathbf{n^1},\mathbf{n^2}),
\label{ntheta}
\end{equation}
where $\mathbf{n^1}$ and $\mathbf{n^2}$ are constrained as in the previous section. Consider the two limiting planes $\mathrm{SFP_1}$ and $\mathrm{SFP_2}$ with normals $\mathbf{n^1}$ and $\mathbf{n^2}$ satisfying the minimum in equation (\ref{ntheta}). Using the rotation argument of Remark \ref{remrot} we can easily show that $\mathrm{SFP_1}$ and $\mathrm{SFP_2}$ are necessary in contact with at least two vertices. However, for each value of the pair $(\theta,\phi)$, we have three types of possible contact between the tetrahedron formed by points $\{{X}^i\}_{i=1\dots 4}$ and the limiting planes: vertex-vertex, edge-vertex, edge-edge or face-vertex.
These configurations can be analyzed as follows: let us consider a pair $(\theta_0,\phi_0)$ and the corresponding type of contact:
\begin{itemize}
\item \textbf{Vertex-vertex}: each limiting plane is in contact with only one vertex, respectively $V_{i}$ and $V_{j}$. In this case, $n(\theta,\phi)$ is differentiable at $(\theta_0,\phi_0)$ and there exists a curve $\gamma_{vv}$ defined by
$$\gamma_{vv}(t)=(\theta(t),\phi(t)),~\gamma_{vv}(0)=(\theta_0,\phi_0),
$$ such that $\frac{d}{dt}n(\theta(t),\phi(t))$ exists and does not vanish for $t=0$. This can be proved using again the rotation argument of Remark \ref{remrot}. Hence, $n(\theta_0,\phi_0)$ cannot be minimal. 
\item \textbf{Edge-vertex}: one of the two planes is in contact with edge $E_{ij}$ and the other one is in contact with vertex $V_k$. In this case there is still a degree of freedom since the plane in contact with $E_{ij}$ can rotate around this edge in either directions while keeping the other plane in contact with $V_k$ only. If the plane in contact with $E_{ij}$ is $\mathrm{SFP_1}$, its normal $\mathbf{n^1}$ can be parameterized by using a single scalar parameter $t$ and we obtain a family of planes defined by
$$
\mathrm{SFP_1}(t)=\left\{\mathbf{X}\in\mathbb{R}^3,\;\scal{\mathbf{n^{1}}(t)}{(\mathbf{X}-\mathbf{X}^i)}=0\right\}.
$$
For each value of $t$, the intersection of $\mathrm{SFP_1}(t)$ with $\mathrm{PSLP}$ defines a Hinge Line and thus a pair $(\theta(t),\phi(t))$ of tilt and swing angles. Hence, there exists a parametric curve 
\begin{equation}\gamma_{ev}(t)=(\theta(t),\phi(t)),~\gamma_{ev}(0)=(\theta_0,\phi_0),
\label{gammaev}
\end{equation}
along which $n(\theta,\phi)$ is differentiable. As we will see in the numerical results, the curve $\gamma_{ev}(t)$ is almost a straight line when $\theta$ and $\phi$ are small, and $\frac{d}{dt}n(\theta,\phi(t))$ does not vanish for $t=0$.
\item\textbf{Edge-edge}: the limiting planes are respectively in contact with edges $E_{ij}$, $E_{kl}$ connecting, respectively, vertices $V_i,V_j$ and vertices $V_k,V_l$. There is no degree of freedom left since these edges cannot be parallel (otherwise all points would be coplanar). Hence, $n(\theta,\phi)$ is not differentiable at ($\theta_0,\phi_0$).
\item \textbf{Face-vertex}: the limiting planes are respectively in contact with vertex $V_l$ and
with the face $F_{ijk}$ connecting vertices $V_i,V_j,V_k$. As in the previous case, there is no degree of freedom left and $n(\theta,\phi)$ is not differentiable at ($\theta_0,\phi_0$).
\end{itemize}

We can already speculate that the first two configurations are necessary suboptimal. Consequently we just have to compute the f-number associated with each one of the 7 possible configurations of type edge-edge of face-vertex.

\subsubsection{Numerical results}
We have considered the flat object of Example \ref{example1_2d}, translated in plane $X_1=-0.5$, and a complimentary point in order to form a tetrahedron. The vertices have coordinates
$$
\mathbf{X}^1=\left(\begin{array}{r}-0.5\\-1\\1\end{array}\right),~\mathbf{X}^2=\left(\begin{array}{r}-0.5\\3\\1\end{array}\right),$$
$$
 \mathbf{X}^3=\left(\begin{array}{r}-0.5\\0\\1.5\end{array}\right),\mathbf{X}^4=\left(\begin{array}{r}1\\1\\1.5\end{array}\right).
$$
\begin{table}
\begin{center}
\begin{tabular}{|c|c|c|c|}
\hline Contact & $\theta$ & $\phi$ & $n(\theta,\phi)$\\
\hline $E_{12}E_{34}$ & 0.021430 &-0.028582& 12.35\\
\hline $E_{23}E_{14}$ & 0.150568 &-0.203694& 90.75\\
\hline $E_{13}E_{24}$ & 0.030005 &-0.020010& 8.64\\
\hline $F_{123}V_{4}$ &  0 & 0.100167 & 88.18\\
\hline $F_{243}V_{1}$ & 0.075070 & -0.050162 &42.91\\
\hline $F_{134}V_{2}$ & 0.033340 & -0.033358 &9.64\\ 
\hline $F_{124}V_{3}$ & 0.018751 &-0.012503 &10.76\\
\hline
\end{tabular}
\end{center}
\caption{Value of $\theta$ and $\phi$ for each possible optimal contact and corresponding f-number $n(\theta,\phi)$.}
\label{table3D}
\end{table}
All configurations of type edge-edge and face-vertex have been considered and 
the corresponding values of $\theta,\phi$ and $n(\theta,\phi)$ are given in Table \ref{table3D}. The $E_{13}E_{24}$ contact seems to give the minimum f-number.
\begin{figure}
\begin{center}
\includegraphics[width=\linewidth]{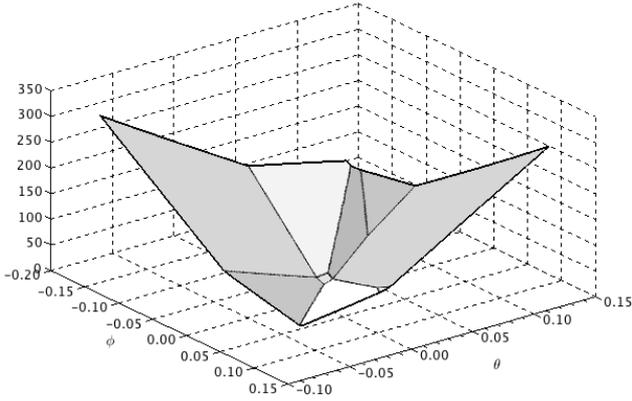}
\end{center}
\caption{Graph of $n(\theta,\phi)$ for $(\theta,\phi)\in[-0.175,0.225]\times[-0.3225,0.2775]$.}
\label{figlevel3d1}
\end{figure}
\begin{figure}
\begin{center}
\includegraphics[width=\linewidth]{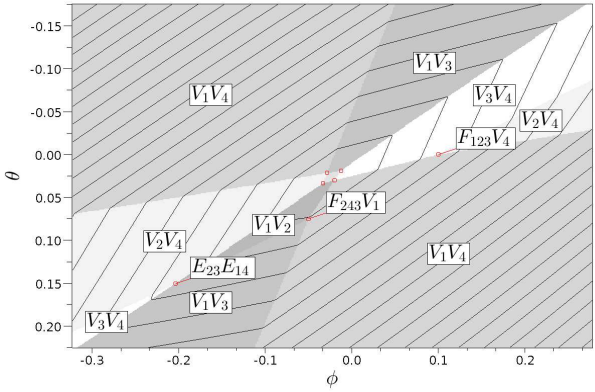}
\end{center}
\caption{Level curves of $n(\theta,\phi)$ and types of contact for $(\theta,\phi)\in[-0.175,0.225]\times[-0.3225,0.2775]$.}
\label{figlevel3d2}
\end{figure}
\begin{figure}
\begin{center}
\includegraphics[width=\linewidth]{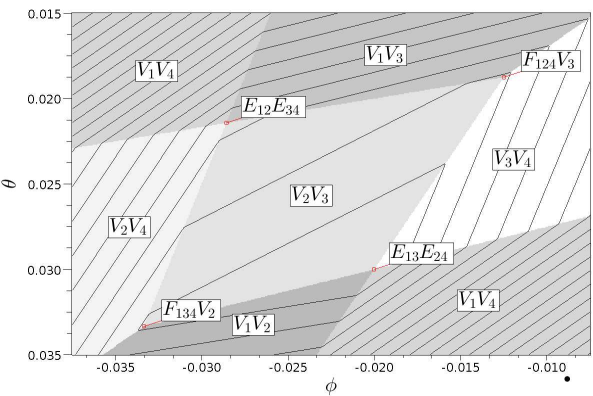}
\end{center}
\caption{Level curves and types of of $n(\theta,\phi)$ contact for $(\theta,\phi)\in[-0.015,0.035]\times[-0.0375,-0.0075]$.}
\label{figlevel3d3}
\end{figure}
Figure \ref{figlevel3d1} gives the graph of $n(\theta,\phi)$ and its level curves in the vicinity of the minimum are depicted in Figures \ref{figlevel3d2} and \ref{figlevel3d3}. In the interior of each different shaded region, the $(\theta,\phi)$ pair is such that the contact of $\mathcal{X}$ with the limiting planes is of type $V_iV_j$. The possible optimal  $(\theta,\phi)$ pairs, corresponding to contacts of type $E_{ij}E_{kl}$ of $F_{ijk}V_l$, are marked with red dots. Notice that the graph of $n(\theta,\phi)$ is almost polyhedral, i.e. in the interior of regions of type $V_iV_j$, the gradient is almost constant and does not vanish, as seen on the level curves. If confirms that the minimum cannot occur in these regions, as announced in Section 3.4.2. 

The frontiers between regions of type $V_iV_j$ are curves corresponding to contacts of type $E_{ij}V_{k}$ and defined by Equation (\ref{gammaev}). The extremities of these curves are $(\theta,\phi)$ pairs corresponding to contacts of type $E_{ij}E_{kl}$ or $F_{ijk}V_l$. For example, in Figure \ref{figlevel3d3}, the $(\theta,\phi)$ pairs on the curve separating $V_2V_3$ and $V_3V_4$ regions correspond to the $E_{24}V_3$ contact. The extremities of this curve are the two $(\theta,\phi)$ pairs corresponding to contacts $F_{124}V_3$ and $E_{13}E_{24}$. Along this curve, $n(\theta,\phi)$ is strictly monotone as shown by its level curves. 

Finally, the convergence of its level curves in Figure \ref{figlevel3d3} confirms that the minimum of $n(\theta,\phi)$ is reached for the $E_{13}E_{24}$ contact. Hence, the optimal angles are $(\theta^*,\phi^*)=(0.030005,-0.020010)$ and the minimal f-number is equal to $n(\theta^*,\phi^*)=8.64$. By comparison, the f-number without tilt and swing optimization is $n(0,0)=28.74$. This example highlights again the important gain in terms of f-number reduction with the optimized tilt and swing angles. In our experience, the optimal configuration for general polyhedrons can be of type edge-edge or face-vertex.

\section{Trends and conclusion}
\label{conclusion}
In this paper, we have given the optimal solution of the most challenging issue in view camera photography: bring an object of arbitrary shape into focus and at the same time minimize the f-number. This problem takes the form of a continuous optimization problem where the objective function (the f-number) and the constraints are non-linear with respect to the design variables. When the object is a convex polyhedron, we have shown that this optimization problem does not need to be solved by classical methods. Under realistic hypotheses, the optimal solution always occurs when the maximum number of constraints are saturated. Such a situation corresponds to a small number of configurations (seven when the object is a tetrahedron). Hence, the exact solution is found by comparing the value of the f-number for each configuration.

The linear algebra framework allowed us to efficiently implement the algorithms in a numerical computer algebra software. The camera software is able to interact with a robotised view camera prototype, which is actually used by our partner photographer. With the robotised camera, the time elapsed in the focusing process is often orders of magnitude smaller than the systematic trial and error technique. 

The client/server architecture of the software allows us to rapidly develop new problem solvers by validating them first on a virtual camera before implementing them on the prototype. We are currently working on the fine calibration of some extrinsic parameters of the camera, in order to improve the precision of the acquisition of 3D points of the object. 

\begin{acknowledgements}
This work has been partly funded by the Innovation and Technological Transfer Center of R\'egion Ile de France.
\end{acknowledgements}

\appendix
\section{Appendix}
\subsection{Computation of the depth of field region}
\label{dofdemo}
\begin{figure}
\begin{center}
\includegraphics[width=\linewidth]{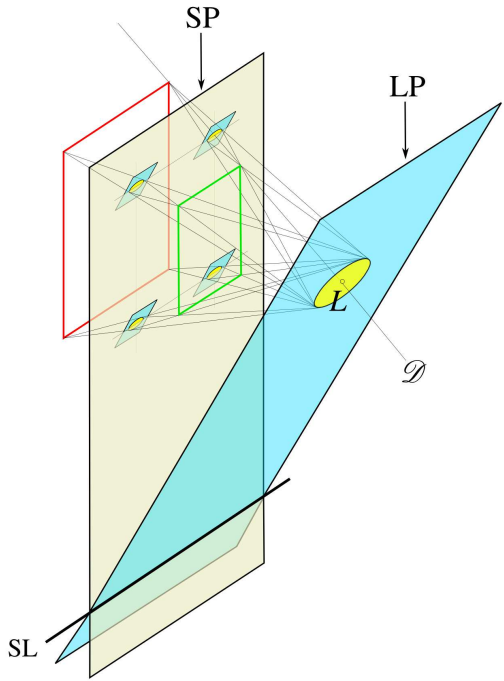}\\
\end{center}
\caption{Construction of four approximate intersections of cones with $\mathrm{SP}$}
\label{figdofdemo}
\end{figure}

\begin{figure}
\begin{center}
\includegraphics[width=\linewidth]{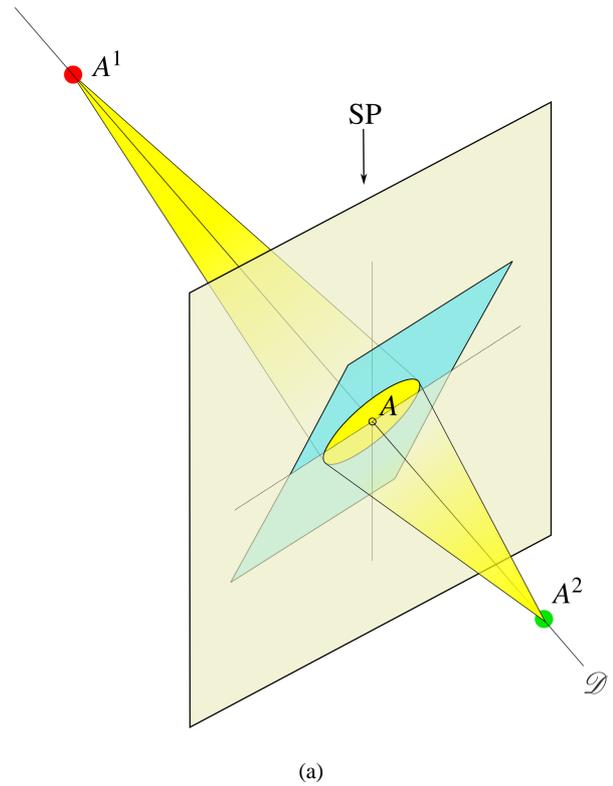}\\
(a)\\
\includegraphics[width=\linewidth]{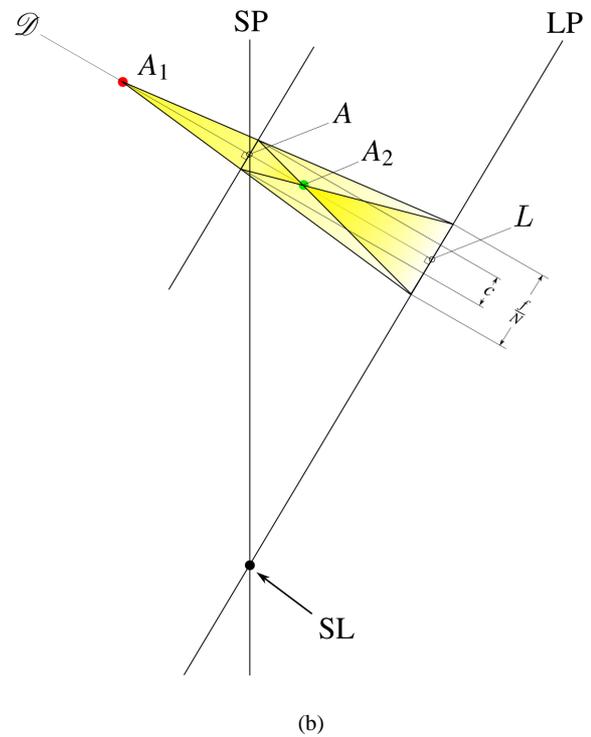}\\
(b)
\end{center}
\caption{(a) close-up of a particular intersection exhibiting the vertices of cones for a given directrix $\mathcal{D}$. (b) geometrical construction allowing to derive the depth of field formula.}
\label{figdofdemo2D}
\end{figure}

In order to explain the kind of approximation used, we have represented in Figure \ref{figdofdemo} the geometric construction of the image space limits corresponding to the depth of field region. Let us consider the cones whose base is the pupil and having an intersection with $\mathrm{SP}$ of diameter $c$. The image space limits are the locus of the vertex of such cones. The key point, suggested in \cite{Bigler} and \cite{Evens}, is the way the diameter of the intersection is measured. 

For a given line $\mathcal{D}$ passing through $L$ and a circle $\mathcal{C}$ of center $L$ in $\mathrm{LP}$
let us call $\mathcal{K}(\mathcal{C},\mathcal{D})$ the set of cones with directrix $\mathcal{D}$ and base 
$\mathcal{C}$. For a given directrix $\mathcal{D}$ let us call $A$ its intersection with $SP$, as depicted in Figure \ref{figdofdemo2D}a. Instead of considering the intersection of cones of directrix $\mathcal{D}$ with $\mathrm{SP}$, we consider their intersections with the plane passing through $A$ and parallel to $\mathrm{LP}$. By construction, all intersections are circles, and there exists only two cones $K_1$ and $K_2$ in $\mathcal{K}(\mathcal{C},\mathcal{D})$ such that this intersection has a diameter equal to $c$, with their respective vertices $A_1,A_2$ on each side of $\mathrm{SP}$, respectively marked in Figure  \ref{figdofdemo2D}a by a red and a green spot. Moreover, for all cones in $\mathcal{K}(\mathcal{C},\mathcal{D})$ only those with a vertex lying on the segment $[A_1,A_2]$ have an "approximate" intersection of diameter less that $c$. 

The classical laws of homothety show that for any directrix $\mathcal{D}$, the locus of the vertices of cones $K_1$ and $K_2$ will be on two parallel planes located in front of and behind $\mathrm{SP}$, as illustrated by a red and a green frame in Figure \ref{figdofdemo}a. Hence, the depth of field region in the object space is the reciprocal image of the region between parallel planes $\mathrm{SP_1}$ and $\mathrm{SP_2}$ as depicted in Figure \ref{scheimpflug-dof}. 

Formulas (\ref{dof}) and (\ref{dofharm}) are obtained by considering the directrix that is orthogonal to $\mathrm{LP}$, as depicted in Figure \ref{figdofdemo2D}b. If we note $p=AL$, $p_1=A_1L$, $p_2=A_2L$, by considering similar triangles, we have
\begin{align}
\frac{p_1-p}{p_1}=\frac{p-p_2}{p_2}=\frac{Nc}{f},
\label{similar}
\end{align}
which gives immediately 
$$
p=\frac{2p_1p_2}{p_1+p_2},
$$
and by substituting $p$ in (\ref{similar}), we obtain 
$$
\frac{p_1-p_2}{p_1+p_2}=\frac{Nc}{f},
$$
which allows to obtain (\ref{dof}).

\subsection{Proof of Proposition \ref{prop1}}
\label{proof:prop1}
Without loss of generality, we consider that the optimal $\theta$ is positive. Suppose now that only one constraint is active in (\ref{cf1}). Then there exists $i_1$ such that $\scal{(\mathbf{X}^{i_1}-\mathbf{W})}{\mathbf{n^1}}=0$ and the first order optimality condition is verified: if we define 
$$
g(\mathbf{a},\theta)=-\scal{(\mathbf{X}^i-\mathbf{W}(\theta))}{\mathbf{n^1}}=X^{i_1}_2-a_1X^{i_1}_3+\frac{f}{\sin\theta},
$$
there exists $\lambda_1\geq 0$ such that the Kuhn and Tucker condition
$$
\nabla N(\mathbf{a},\theta)+\lambda_1\nabla g(\mathbf{a},\theta)=0,
$$
is verified. Hence, we have
$$
\frac{2}{(2\cot \theta -(a_1+a_2))^2}\left(\begin{array}{c}\cot\theta-a_2\\-\cot\theta+a_1\\\frac{a_1-a_2}{\sin^2\theta}\end{array}\right)+\lambda_1\left(\begin{array}{c}-X^{i_1}_3\\0\\-\frac{\cos\theta}{\sin^2\theta}\end{array}\right),
$$
and necessarily, $a_1=\cot\theta$ so that $N(a,\theta)$ reaches its upper bound and thus is not minimal. We obtain the same contradiction when only a constraint in (\ref{cf2}) is active, or only two constraints in (\ref{cf1}), or only two constraints in (\ref{cf2}).\qed

\subsection{Proof of Proposition \ref{prop2}}
\label{proof:prop2}
Without loss of generality we suppose that $\theta\geq 0$. Suppose that the minimum of $N(\mathbf{a},\theta)$ is reached with only vertices $i_1$ and $i_2$ respectively in contact with limiting planes $\mathrm{SFP_1}$ and $\mathrm{SFP_2}$.  The values of $a_1$ and $a_2$ can be determined as the following functions of $\theta$
$$
a_1(\theta)=\frac{X_2^{i_1}+\frac{f}{\sin\theta}}{X_3^{i_1}},~a_1(\theta)=\frac{X_2^{i_2}+\frac{f}{\sin\theta}}{X_3^{i_2}},
$$
and straightforward computations give
\begin{equation}
N(\mathbf{a}(\theta),\theta)=\frac{\left(\frac{X^{i_1}_2}{X^{i_1}_3}-\frac{X^{i_2}_2}{X^{i_2}_3}\right)\sin\theta +\left(\frac{f}{X^{i_1}_3}-\frac{f}{X^{i_2}_3}\right)}{2\cos\theta-\left(\frac{X^{i_1}_2}{X^{i_1}_3}+\frac{X^{i_2}_2}{X^{i_2}_3}\right)\sin\theta-\left(\frac{f}{X^{i_1}_3}+\frac{f}{X^{i_2}_3}\right)}\left(\frac{f}{c}\right).
\label{dofproof}
\end{equation}
In order to prove the result, we just have to check that the derivative of $N(\mathbf{a}(\theta),\theta)$ with respect to $\theta$ cannot vanish for $\theta\in[0,\frac{\pi}{2}]$. The total derivative of $N(\mathbf{a}(\theta),\theta))$ with respect to $\theta$ is given by

\begin{multline*}\frac{d}{d\theta}N(\mathbf{a}(\theta),\theta)=\\
2\frac{\left(\frac{X^{i_1}_2}{X^{i_1}_3}-\frac{X^{i_2}_2}{X^{i_2}_3}\right)-f\left(\frac{X^{i_1}_2-X^{i_2}_2}{X^{i_1}_3 X^{i_2}_3}\right)\cos\theta+\left(\frac{f}{X^{i_1}_3}-\frac{f}{X^{i_2}_3}\right)\sin\theta}{\left(2\cos\theta-\left(\frac{X^{i_1}_2}{X^{i_1}_3}+\frac{X^{i_2}_2}{X^{i_2}_3}\right)\sin\theta-\left(\frac{f}{X^{i_1}_3}+\frac{f}{X^{i_2}_3}\right)\right)^2}\left(\frac{f}{c}\right)
\end{multline*}
and its numerator is proportional to the trigonometrical polynomial
\begin{equation}
p(\theta)=b_0+b_1\cos\theta+b_2\sin\theta, 
\label{ptheta}
\end{equation}
where $b_0=X^{i_1}_2X^{i_2}_3-X^{i_2}_2X^{i_1}_3$, $b_1=-f\left(X^{i_2}_2-X^{i_1}_2\right)$, $b_2=f(X^{i_2}_3-X^{i_1}_3)$. It can be easily shown by using the Schwartz inequality that $p(\theta)$ does not vanish provided that 
\begin{equation}
b_0^2>b_1^2+b_2^2.
\label{eqabc}
\end{equation}
Since $X_1^{i_2}=X_1^{i_1}=0$, whe have $b_0^2=\Vert X_1^{i_1}\times X_1^{i_1}\Vert^2$ and $b_1^2+b_2^2=f^2\Vert X_1^{i_1}- X_1^{i_1}\Vert^2$. Hence (\ref{eqabc}) is equivalent to condition (\ref{conddet}), this ends the proof.\qed

\subsection{Depth of field region approximation used by A. Merklinger}
\label{appendix1}

\begin{figure}
\begin{center}
\includegraphics[width=\linewidth]{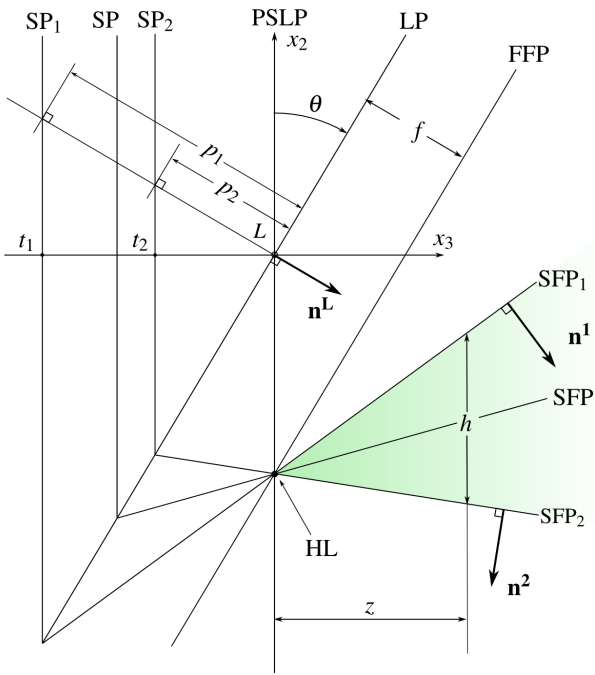}
\end{center}
\caption{Depth of field region approximation using the hyperfocal.}
\label{dofmerk}
\end{figure}
In his book (\cite{Merklinger}, Chapter 7) A. Merklinger has proposed the following approximation based on the assumption of distant objects and small tilt angles: if $h$ is the distance from $\mathrm{SFP_1}$ to $\mathrm{SFP_2}$, measured in a direction parallel to the sensor plane at distance $z$ from the lens plane, as depicted in Figure \ref{dofmerk}, we have
\begin{equation}
\frac{h}{2z}\approx\frac{f}{H\sin\theta},
\label{dofhyp}
\end{equation}
where $H$ is the \textit{hyperfocal distance} (for a definition see \cite{Ray} p. 221), related to the f-number $N$ by the formula
$$
H=\frac{f^2}{Nc},
$$
and $c$ is the diameter of the circle of confusion. Using (\ref{dofhyp}), the f-number $N$ can be approximated by
$$
\tilde N=\sin\theta\frac{h}{2z}\frac{f}{c}.
$$
Since the slopes of $\mathrm{SFP_1}$ and $\mathrm{SFP_2}$ are respectively given by $a_1$ and $a_2$, we have 
$$
\frac{h}{z}={a_1-a_2},
$$
and we obtain immediately
$$
\tilde N=(a_1-a_2)\sin\theta\left(\frac{f}{2c}\right),
$$
which is the same as (\ref{approxN}).

\bibliographystyle{ieee}
\bibliography{egbib}

\end{document}